\documentclass[12pt]{article}

\usepackage{shortcuts}
\usepackage{url}
\usepackage{graphicx}
\usepackage{hyperref}
\usepackage{multirow}
\usepackage{subcaption}

\hypersetup{
	colorlinks,
	citecolor=black,
	filecolor=black,
	linkcolor=black,
	urlcolor=black
}

\newcommand{\quand}{\quad \text{ and } \quad}

\newcommand{\tZ}{{\tilde{Z}}}

\begin{document}
\baselineskip=2pc
\begin{center}
	\Large{\bf Enforcing strong stability of explicit Runge--Kutta methods with superviscosity}
\end{center}
\centerline{
	Zheng Sun\footnote{Department of Mathematics, The Ohio State University,
		Columbus, OH 43210, USA. E-mail: sun.2516@osu.edu.} and 
	Chi-Wang Shu\footnote{Division of Applied Mathematics, Brown University,
		Providence, RI 02912, USA. E-mail: chi-wang\_shu@brown.edu.
Research supported by NSF grant DMS-1719410.}
}

\centerline{\bf Abstract} 
\smallskip

A time discretization method is called strongly stable, if the norm of its numerical solution is nonincreasing. It is known that, even for linear semi-negative problems, many explicit Runge--Kutta (RK) methods fail to preserve this property. In this paper, we enforce strong stability by modifying the method with superviscosity, which is a numerical technique commonly used in spectral methods. We propose two approaches, the modified method and the filtering method for stabilization. The modified method is achieved by modifying the semi-negative operator with a high order superviscosity term; the filtering method is to post-process the solution by solving a diffusive or dispersive problem with small superviscosity. For linear problems, most explicit RK methods can be stabilized with either approach without accuracy degeneration. Furthermore, we prove a sharp bound (up to an equal sign) on diffusive superviscosity for ensuring strong stability. The bound we derived for general dispersive-diffusive superviscosity is also verified to be sharp numerically. For nonlinear problems, a filtering method is investigated for stabilization. Numerical examples with linear non-normal ordinary differential equation systems and for discontinuous Galerkin approximation of conservation laws are performed to validate our analysis and to test the performance. 

\vspace{.2in}

 \bigskip

 \bigskip

\clearpage
\section{Introduction}
\setcounter{equation}{0}
\setcounter{figure}{0}
\setcounter{table}{0}

When discretizing hyperbolic conservation laws, a stable spatial discretization may result in a method-of-lines scheme of the form
\begin{equation}\label{eq-nlode}
\frac{d}{dt}u(t) = F(u(t)), \qquad  \ip{F (v), v}  \leq 0, \qquad \forall v\in V,
\end{equation}
where $(V,\ip{\cdot,\cdot})$ is a real inner-product space. Under the associated norm $\|\cdot\| = \sqrt{\ip{\cdot,\cdot}}$, the solution to \eqref{eq-nlode} obeys the energy decay law
\begin{equation}\label{eq-decay}
\frac{d}{dt}\nm{u}^2 = 2\ip{F(u),u} \leq 0.
\end{equation}
For time integration of \eqref{eq-nlode} with Runge--Kutta (RK) methods, if the monotonicity of the solution norm is preserved in the discrete sense
\begin{equation}
\nm{u^{n+1}}\leq \nm{u^n},
\end{equation}
possibly under an appropriate time step constraint, then the method is called \emph{strongly stable}. It is also referred to as \emph{monotonicity-preserving} or \emph{contractive} in the ordinary differential equation (ODE) literature.

Compared with a weaker type of stability 
$\nm{u^{n+1}}\leq (1+C\tau)\nm{u^n}$ (where $\tau$ is the time step size), 
strong stability has received particular attention in certain scenarios. Firstly, the weaker stability implies $\nm{u^n}\leq Ce^{CT} \nm{u^0}$, which means the solution norm may grow exponentially with respect to the final time $T$ \cite{gustafsson1995time}. The strong stability, which offers $\nm{u^n}\leq \nm{u^0}$,  maintains better control on the solution after long time integration. Secondly, strong stability also connects to contractivity of the time integrator \cite{butcher2016numerical}. It guarantees that the numerical error in each time step will not be amplified during the time integration. Finally, for discretizing conservation laws, the strong stability retrieves the entropy inequality under the square entropy functional. Obeying such entropy inequality is crucial for capturing the physically relevant solution. Numerous efforts have been dedicated to semi-discrete entropy stable schemes \cite{tadmor1987numerical, fjordholm2012arbitrarily, fisher2013high, chen2017entropy, abgrall2018general, chen2019review}. However, this property may not be preserved after time discretization \cite{lozano2018entropy,lozano2019entropy}. Enforcing strong stability is usually the first step towards a fully discrete entropy stable scheme.

For hyperbolic conservation laws and many other problems, explicit RK methods are widely used due to its simplicity and capability to be incorporated with limiters for suppressing spurious oscillations and preserving physical bounds. See \cite{gottlieb2001strong, cockburn2001runge,zhang2010maximum,zhang2010positivity,xing2010positivity,wu2015high,qin2016bound,guo2017bound,sun2018discontinuous,sun2019entropy} and references therein. Compared with implicit time integrators (such as algebraically stable methods \cite{butcher2016numerical}), strong stability of explicit RK methods is not well understood until recently. Classical eigenvalue analysis can be used for analyzing linear problems with normal operators, but they may fail for non-normal problems. It could give a wrong bound for time step constraints \cite{iserles2009first}, or even misjudge the strong stability of a method \cite{sun2017rk4}.  The energy methods are hence used for analyzing non-normal problems. For coercive operators, after the earlier work by Levy and Tadmor \cite{levy1998semidiscrete}, in \cite{gottlieb2001strong}, by exploiting the strong-stability-preserving formulation, the authors showed that all linear RK methods are strongly stable. This coincides with the analysis on conditional contractivity of explicit RK methods in the ODE literature \cite{spijker1983contractivity, kraaijevanger1991contractivity}. For general linear noncoercive operators, besides analysis for particular RK schemes \cite{tadmor2002semidiscrete,sun2017rk4,ranocha2018l_2}, in \cite{sun2019strong}, a unified framework is established to analyze RK methods of arbitrary order. In particular, we proved that for $p$-stage $p$th order linear RK methods, the methods are not strongly stable if $p\equiv 1 \Mod 4$ or $p\equiv 2 \Mod 4$, and the methods are strongly stable if $p\equiv 3 \Mod 4$. Another stream of studies particularly focuses on discontinuous Galerkin (DG) operators for linear conservation laws. Besides stability results that are consistent with those for general operators, it is also shown that the methods can achieve strong stability with extra low order spatial polynomials  \cite{zhang2010stability,sun2017stability, xu2019L2,xu2020superconvergence}. For nonlinear problems, counterexamples were given by Ranocha in \cite{ranocha2018strong} showing that various 
explicit RK methods do not preserve strong stability and order barriers were proved for certain strong-stability-preserving RK methods. Further recent studies on criteria for nonlinear problems can be found in \cite{ranocha2019energy}.

Due to the importance of strong stability and the negative fact that the property is not typically preserved, it is natural to pursue remedy approaches to stabilize explicit RK methods. In \cite{ketcheson2019relaxation}, Ketcheson proposed the relaxation RK methods to enforce strong stability, which introduces a relaxation factor computed from intermediate stages for updating the solution increment. Soon in \cite{ranocha2019relaxation}, this method is generalized to convex functionals for designing fully discrete entropy stable methods for Euler and Navier-Stokes equations. In this paper, we investigate another stabilization approach based on \emph{superviscosity} (or  \emph{hyperviscosity}), a numerical  technique commonly used in spectral methods. See \cite{passot1988hyperviscosity,tadmor1993super, guermond2003mathematical, tadmor2004burgers} and references therein. The idea for superviscosity is to add an additional high order derivative term, which penalizes the high-frequency modes of the solution to provide extra regularization and to damp the solution norm. For scalar conservation laws, the modified equation  
\begin{equation}\label{eq-supvisc}
\frac{\partial }{\partial t} u +\frac{\partial }{\partial x} f(u) = \veps_\nu (-1)^{\kappa+1}  \frac{\partial^{2\kappa} }{\partial x^{2\kappa}} u, \qquad 0<\veps_\nu \ll 1
\end{equation}
is solved instead of the original inviscid equation. 
More generally, a convolution kernel can be inserted into the superviscosity term, namely $ \veps_\nu  \left(-\frac{\partial }{\partial x}\right)^\kappa \left(Q_m(x,t)* \left(\frac{\partial }{\partial x}\right)^\kappa u\right)$, to tune the inviscid spectrum and the penalization on high frequency modes. This approach is known as spectral (super-/hyper-)viscosity method and the case $\kappa = 1$ has received particular interest  \cite{tadmor1989convergence, maday1989analysis, maday1993legendre,chen1993spectral,karamanos2000spectral,tadmor2012adaptive}. Furthermore, it is known that the superviscosity approach is closely related with the low pass filter \cite{gottlieb2001spectral, hesthaven2007spectral}, which post-processes the solution after each time step (or every a few time steps) by directly modifying Fourier coefficients. 

In this paper, we generalize the idea of superviscosity in spectral methods and incorporate it with time integrators. Our main focus is on linear autonomous problems with non-normal operators. We consider not only the diffusive superviscosity in \eqref{eq-supvisc}, but also the dispersive superviscosity with odd order derivatives. The dispersive superviscosity is shown to be necessary to stabilize certain RK methods (for example, the fourth order RK method). Two approaches are investigated in parallel. The first one is the modified method, which supplements the original semi-negative operator with a superviscosity term and then use the classical RK methods for time integration. The second approach is the filtering method, which post-processes the solution of the original RK methods by solving a diffusive and (or) dispersive problem with small superviscosity. Both methods are shown to be able to stabilize RK methods for linear problems. Furthermore, we provide a bound on superviscosity that is needed to enforce strong stability. Up to an equal sign, the bound for diffusive superviscosity is proved to be sharp and the sharpness for the general dispersive-diffusive case is also confirmed numerically. Moreover, the inclusion of superviscosity in the time integrator will not cause accuracy degeneration. For nonlinear problems, a filtering approach with diffusive superviscosity is considered. The amplitude of the superviscosity term is calculated from the previous solution and the preliminary solution. Note the discretization in this paper would retrieve (or be very close to) the classical approaches with spectral methods for spatial approximations, whose performance has been well studied in the literature. To avoid repetition, we will use an ODE problem and DG approximations of conservation laws for numerical examples. 

The rest of the paper is organized as follows. In Section \ref{sec-pre}, we review estimates in \cite{sun2019strong} for understanding the energy change of classical explicit RK methods in each time step.  In Section \ref{sec-mf}, the modified and the filtering methods are presented for linear problems. We estimate the energy change between the solution given by classical RK methods and the stabilized solution. Then we combine it with results in Section \ref{sec-pre} to conclude the superviscosity approach is able to enforce strong stability without affecting accuracy order. The technical proof of Lemma \ref{lem-key} is given in Section \ref{sec-thms}. In Section \ref{sec-nonlinear}, we briefly discuss a filtering method for nonlinear problems. Numerical tests are presented in Section \ref{sec-num} to validate our analysis and the conclusion is given in Section \ref{sec-concl}. 

\section{Preliminaries: RK methods for linear ODEs}\label{sec-pre}
\setcounter{equation}{0}
\setcounter{figure}{0}
\setcounter{table}{0}

In this section, we review some of the results in \cite{sun2019strong}, which sets the foundation for analysis in this paper. 

Suppose $F(u) = Lu$ is linear in \eqref{eq-nlode}. In other words, we consider the semi-negative  linear autonomous system 
\begin{equation}\label{eq-ode}
\frac{d}{dt}u = Lu, \qquad  \ip{L v, v}  \leq 0, \qquad \forall v \in V. 
\end{equation}
A consistent explicit RK method for \eqref{eq-ode} with the input $u_0$ and time step $\tau>0$ takes the form
\begin{equation}\label{eq-rk}
u_+ = R_s(Z) u_0, \qquad R_s(Z) = \sum_{k = 0}^s \alpha_k Z^k,\quad Z = \tau L, \quand \alpha_0 = \alpha_1 = 1.
\end{equation}
The method is (linearly) $p$th order if and only if $\alpha_k = \frac{1}{k!}$, $\forall k\leq p$ and $\alpha_{p+1}\neq \frac{1}{(p+1)!}$. 

By taking inner-product with itself on both sides and subtracting 
$\nm{u_0}^2$, one can see the energy change after one time step of \eqref{eq-rk} is 
\begin{equation}\label{eq-sq}
	\nm{u_+}^2 - \nm{u_0}^2 = \sum_{i,j = 0,~ i^2+j^2 \neq 0}^s \alpha_i \alpha_j \ip{Z^i u_0, Z^j u_0}.
\end{equation}
We introduce the semi-inner-product\footnote{Note in \cite{sun2019strong}, we exclude the time step $\tau$ in the definition of semi-inner-product $\qip{\cdot,\cdot}$. In other words, $L$ instead of $Z$ is used to define $\qip{\cdot,\cdot}$. }
\begin{equation}
\qip{v,w}:= - \ip{Zv,w} - \ip{v,Zw}
\end{equation}
and denote by $\qnm{v} = \sqrt{\qip{v,v}}$ the induced semi-norm.
After repeatedly using the identity
\begin{equation}\label{eq-dibp}
	\ip{Zv,w} = -\ip{v,Zw} - \qip{v,w},
\end{equation}
one can rewrite the energy equality \eqref{eq-sq} into the form
\begin{equation}\label{eq-ener}
\nm{u_+}^2 - \nm{u_0}^2 =  \sum_{k = 1}^s\beta_k  \nm{Z^k u_0}^2 + \sum_{i,j = 0}^{s-1}\gamma_{i,j} \qip{Z^i u_0, Z^j u_0}, \qquad \gamma_{i,j} = \gamma_{j,i}.
\end{equation}
\eqref{eq-dibp} is the discrete analogue of integration by parts when $L$ approximates $\partial_x$. The conversion from \eqref{eq-sq} to \eqref{eq-ener} can be achieved by using a simple computer program. See \cite[Algorithm 1]{sun2019strong} for details.

As has been shown in \cite{sun2019strong}, for $\|Z\|=\tau \|L\| \leq \lambda \ll 1$, the sign of the energy change is determined by its low order terms. We introduce the definition below to facilitate our discussion. 

\begin{DEFN}
	The \emph{leading index}, $k^*$, is the index of the first nonzero $\beta_k$. In other words, $\beta_{k^*} \neq 0$ and $\beta_k = 0$, $\forall 1\leq k < k^*$. The \emph{leading submatrix} $\Gamma^* := \Gamma_{k^*-1} = (\gamma_{i,j})_{0\leq i,j\leq k^*-1}$ is the $k^*$th order symmetric matrix given 
by the coefficients $\gamma_{i,j}$.
\end{DEFN} 

Throughout the paper, we use $\xi(\cdot)$ and $\tilde{\xi}(\cdot)$ to represent polynomials, whose coefficients may depend on $\{\alpha_k\}$, $\mu$ and $\nu$, but are independent of $\lambda$, $\tau$ and the underdetermined parameter $\veps$. Following lines in \cite{sun2019strong}, one can obtain an energy estimate from \eqref{eq-ener}. 
\begin{PROP}\label{prop-ener} Suppose $E^* = \diag(e_0,\ldots,e_{k^*-1})$ and $\Gamma^* - {E}^*<0$. There exists a constant $\lambda \ll 1$, such that
	\begin{equation}\label{eq-ener-prop}
	\nm{u_+}^2 - \nm{u_0}^2 \leq \left(\beta_{k^*} + \lambda \xi(\lambda)\right)\nm{Z^{k^*} u_0}^2  +\sum_{k = 0}^{k^*-1}e_k\qnm{Z^k u_0}^2,
	\end{equation} 
if $\nm{Z}\leq \lambda$.
\end{PROP}

To better characterize the energy change, we need to know values of $k^*$, $\Gamma^*$ and $\beta_{k^*}$. This has been investigated in  \cite[Lemma 4.1]{sun2019strong}. Proposition \ref{prop-kstar} is its corollary with results that are needed in this paper.

\begin{PROP}\label{prop-kstar}
	\
	\begin{enumerate}
		\item 	If $p$ is odd, then $k^* = \frac{p+1}{2}$ and $\Gamma^* < 0$. 
		\item 	If $p$ is even, then $k^* \geq \frac{p}{2}+1$ and $\Gamma_{\frac{p}{2}-1} = (\gamma_{i,j})_{0\leq i,j\leq \frac{p}{2}-1} < 0$. 
		\item  	For $p$-stage $p$th order linear RK methods, we have 
		\begin{enumerate}
			\item if $p$ is  odd, then $\beta_{k^*} = \frac{(-1)^{ \frac{p+1}{2}}}{(p+1)!}$; 
			\item if $p$ is  even, then $k^* = \frac{p}{2} + 1$ and $\beta_{k^*} = \frac{(-1)^{ \frac{p}{2}+1}}{p!(p+2)}$. 
		\end{enumerate}
	\end{enumerate}

\end{PROP}

\begin{REM}\label{rem-E}
Proposition \ref{prop-kstar} indicates the following choices of $E^*$ in Proposition \ref{prop-ener} for a simplified energy estimate. If $p$ is odd, $\Gamma^*<0$ and Proposition \ref{prop-ener} holds with $E^* = 0$. When $p$ is even and $k^* = \frac{p}{2} + 1$, Proposition \ref{prop-kstar} indicates that there exists $\mu_0$ such that  $\Gamma^* - \diag(0,\cdots,0, \mu) < 0 $ if $\mu>\mu_0$. We can take $E^* =  \diag(0,\cdots,0, \mu)$.
\end{REM}

\section{Modified and filtering methods for linear problems}\label{sec-mf}
\setcounter{equation}{0}
\setcounter{figure}{0}
\setcounter{table}{0}

\subsection{Methods with superviscosity} 

\textbf{The modified method:} For the modified method, we use the original RK method to solve the modified ODE,
\begin{equation}\label{eq-mode}
\dt u = \tL u,
\end{equation}
where
\begin{equation}
\tL  =  L +  \mu \tau^{2k^*-2} (L^\top)^{k^*-1} L^{k^*} + \nu \tau^{2k^*-1} (L^\top)^{k^*} L^{k^*},
\end{equation}
and $L^\top$ is the adjoint operator of $L$ in $V$. 
In other words, we set
\begin{equation}
u_{\rm{M}} = R_s(\tZ) u_0, \qquad \tZ = \tau \tL =  Z +  \mu (Z^\top)^{k^*-1} Z^{k^*} + \nu  (Z^\top)^{k^*} Z^{k^*},
\end{equation}
as the solution at the next time step. The $\mu $ term corresponds to 
dispersive superviscosity and the $\nu$ term corresponds to diffusive superviscosity. 

\begin{REM}
	Note that 
	\begin{equation}
		\ip{\tL v, v} = \ip{Lv,v}- \frac{\mu }{2}\qnm{Z^{k^*-1} u}^2 + \nu \nm{Z^{k^*}v}^2\leq - \frac{\mu }{2}\qnm{Z^{k^*-1} u}^2 + \nu \nm{Z^{k^*}v}^2,~~\forall v\in V
	\end{equation} When $\mu > 0 $ and $\nu <0$, the modified operator is weakly coercive. The superviscosity term introduces extra dissipation in the operator, which helps stabilize the method.
\end{REM}
\textbf{The filtering method:} For the filtering method, we use $u_+$ defined in \eqref{eq-rk} as a preliminary solution, and then apply the superviscosity filter. 
\begin{equation}\label{eq-filter}
	u_{\rm{F}} = \left(I +  \mu (Z^\top)^{k^*-1} Z^{k^*} + \nu  (Z^\top)^{k^*} Z^{k^*}\right) u_+.
\end{equation}

\begin{REM}
The filtering method can be interpreted as a splitting approximation to   \eqref{eq-mode}. We first use the original RK method to solve 
$\frac{d}{dt}u = Lu$ and then use Euler forward method to solve $\frac{d}{dt} u = \left(\tL - L\right)u$.
\end{REM}

It can be seen that the difference between the modified and the filtering method is only a high order term. Hence the performance of two methods is quite similar, especially when the time step is small. This is confirmed by our analysis and numerical tests. From a practical point of view, the filtering method seems to be more efficient.

\subsection{Main results}
In Lemma \ref{lem-key}, we establish an inequality between solution norms of the original method and the stabilized methods. The proof will be given in Section \ref{sec-thms}. 

\begin{LEM}\label{lem-key}
	 For any $\mu$ and $\nu$ with $\mu = 0$ for $p = 1$, there exists a constant $\lambda$, such that for any $0<\veps<1$, 	if $\nm{Z}\leq \lambda$, we have
	\begin{equation}\label{eq-leq}
	 \nm{u_\#}^2- \nm{u_+}^2\leq ( -\mu + \veps |\mu|)\qnm{Z^{{k^*}-1}u_0}^2 + (2\nu +\frac{\lambda}{\veps} \xi(\lambda))\nm{Z^{{k^*}}u_0}^2,
	\end{equation}
	and 
	\begin{equation}\label{eq-geq}
	\nm{u_\#}^2- \nm{u_+}^2\geq (-\mu-\veps |\mu| )\qnm{Z^{{k^*}-1}u_0}^2 + (2\nu +\frac{\lambda}{\veps} \tilde{\xi}(\lambda))\nm{Z^{{k^*}}u_0}^2,
	\end{equation}
	where $ \# \in \{ \rm{F},   \rm{M}\}$.
\end{LEM}

Theorem \ref{thm-key} confirms that, for linear autonomous problems, it is always possible to enforce strong stability by adding sufficiently strong superviscosity. Furthermore, this procedure does not affect accuracy order of the original method. 

\begin{THM}\label{thm-key}
	For an RK method with $k^* = \lceil \frac{p+1}{2}\rceil$, the corresponding modified and filtering methods are (linearly) $p$th order accurate and are strongly stable, if $\nu <\nu_0 = -\frac{\beta_{k^*}}{2}$ and 
\begin{enumerate}
	\item $\mu = 0$ for odd  $p$,
	\item $\mu >\mu_0$ for even $p$.
\end{enumerate}
Here $\mu_0$ is a real number so that $\Gamma^* - \diag(0,\cdots,0,\mu_0) \leq 0$, whose existence has been confirmed in Remark \ref{rem-E}. 
\end{THM}

\begin{proof}
	We first prove strong stability. If $p$ is odd, for $\mu = 0$ and $\lambda \ll1$, we combine \eqref{eq-ener-prop} and \eqref{eq-leq} (with $\veps = 1$ and $E^* = 0$) to get 
	\begin{equation}
	\nm{u_\#}^2 - \nm{u_0}^2 \leq \left(\beta_{k^*}+2\nu + {\lambda} \xi(\lambda)\right)\nm{Z^{{k^*}}u_0}^2,\quad\# \in \{ \rm{F},   \rm{M}\},
	\end{equation}
	The right hand side is negative when $\nu < -\frac{\beta_{k^*}}{2}$ and $\lambda \ll 1$. 
	If $p$ is even, then  we have $k^* = \frac{p}{2}+1$ under our assumption and $\Gamma^* - E^*(\mu_0+\veps)< 0$, $\forall \veps >0$ according to Remark \ref{rem-E}. Combining Proposition \ref{prop-ener} and Lemma \ref{lem-key} gives 
	\begin{equation}
	\nm{u_\#}^2 - \nm{u_0}^2 \leq \left(\veps|\mu|-\mu + \mu_0 + \veps\right)\qnm{Z^{{k^*}-1}u_0}^2 + \left(\beta_{k^*}+2\nu +\frac{\lambda}{\veps} \xi(\lambda)\right)\nm{Z^{{k^*}}u_0}^2,\quad\# \in \{ \rm{F},   \rm{M}\},
	\end{equation}
	if $\nm{Z}\leq \lambda \ll  1$. 
	Therefore, when $\mu> \mu_0$ and $\nu < \nu_0 =- \frac{\beta_{k^*}}{2}$, we can choose $\veps>0$ sufficiently small, such that $\veps|\mu|-\mu + \mu_0 + \veps < 0$, and then take $\lambda\ll 1$ to make the second term negative. Hence the method becomes strongly stable under our choices of $\mu$ and $\nu$. 
	
	To show they are still $p$th order accurate, we note that
	\begin{equation}
	\nm{\mu (Z^\top)^{k^*-1} Z^{k^*} + \nu  (Z^\top)^{k^*} Z^{k^*}} =  \mathcal{O}\left(\tau^{{2k^*-1}}\left(\mu+\nu\tau\right)\right).
	\end{equation}
	Under the assumption $k^* = \lceil \frac{p+1}{2}\rceil$, for odd $p$ with $\mu = 0$ and even $p$ with $\mu = \mathcal{O}(1)$, the superviscosity adds an extra $\mathcal{O}(\tau^{p+1})$ term in the truncation error. Therefore, it will not affect the accuracy order of the time integrator, which completes the proof.
\end{proof}

Theorem \ref{thm-sharp} states the bound for diffusive superviscosity is sharp.
\begin{THM}\label{thm-sharp}
	The modified and filtering methods are NOT strongly stable if $\nu> \nu_0 = -\frac{\beta_{k^*}}{2}$. 
\end{THM}
\begin{proof}
	Consider an energy conservative problem with $\ip{Lv,v} = 0$, $\forall v \in V$. Let $\veps = 1$ in \eqref{eq-geq}, we have
	\begin{equation}
	\nm{u_\#}^2- \nm{u_+}^2 \geq  (2\nu + {\lambda}{\xi}(\lambda))\nm{Z^{{k^*}}u_0}^2.
	\end{equation}
	Then by substituting \eqref{eq-ener} into the estimate above and noting that $\qip{\cdot,\cdot } = 0$, one can see that 
	\begin{equation}
	\begin{aligned}
	\nm{u_\#}^2- \nm{u_0}^2 \geq &\; (\beta_{k^*} + 2\nu + {\lambda} {\xi}(\lambda))\nm{Z^{{k^*}}u_0}^2 + \sum_{k = k^* + 1}^s \beta_k\nm{Z^k u_0}^2 \\
	\geq &\; (\beta_{k^*} + 2\nu + {\lambda} {\xi}(\lambda))\nm{Z^{{k^*}}u_0}^2 - \left(\sum_{k = k^* + 1}^s| \beta_k|\lambda^{k-k^*}\right)\nm{Z^{k^*} u_0}^2 \\
	\geq &\; (\beta_{k^*} + 2\nu + {\lambda} {\xi}(\lambda))\nm{Z^{{k^*}}u_0}^2. 
	\end{aligned}
	\end{equation}
	As $\lambda \ll 1$, and $\nu > - \frac{\beta_{k^*}}{2}$, the right hand side is positive whenever $\nm{Z^{k^*} u_0}\neq 0$. Hence the methods are not strongly stable. 
\end{proof}
\begin{REM}
	In Theorem \ref{thm-key}, whether the equal sign can be taken may depend on the original RK method. The sharpness of $\mu_0$ is verified numerically in Section \ref{sec-num}.
\end{REM}

\section{Proof of Lemma \ref{lem-key}}
\label{sec-thms}
\setcounter{equation}{0}
\setcounter{figure}{0}
\setcounter{table}{0}

In this section, we prove Lemma \ref{lem-key}. Only \eqref{eq-leq} will be considered and \eqref{eq-geq} can be shown following the same blueprint. The main idea is to estimate the high-order terms and keep only the crucial parts in 
\begin{equation}
	\nm{u_\#}^2 = \nm{u_+}^2 + 2\ip{u_+, u_\# - u_+} + \nm{u_\#-u_+}^2.
\end{equation}
We start with introducing a few inequalities that will be used in the proof. Then energy estimate for the filtering method will be proved. Finally, we investigate the modified method.
\subsection{A few inequalities}

\begin{LEM}\label{lem-Zjump}
		\begin{equation}\label{eq-Zjump}
	\qnm{Z^lv} \leq \sqrt{2}\lambda^{l+\hf}\nm{v}.
	\end{equation}
\end{LEM}
\begin{proof}	$
	\qnm{Z^lv}^2 = - \ip{Z^lv,Z^{l+1}v} - \ip{Z^{l+1}v,Z^lv} \leq 2\nm{Z}^{2l+1}\nm{v}^2 \leq 2\lambda^{2l+1}\nm{v}^2$. \eqref{eq-Zjump} can be shown after taking the square root on both sides of the inequality. 
\end{proof}
Lemma \ref{lem-filter} will be used in the proof of the filtering method.
\begin{LEM}\label{lem-filter}
	Let $v_+ = R_s(Z)v$. Then we have 
	\begin{equation}
		\left|\nm{Z^l v_+}^2 -  \nm{Z^l v}^2\right| \leq \lambda\xi(\lambda)\nm{Z^l v}^2,
	\end{equation}
	\begin{equation}
		\left| \qnm{Z^l v_+}^2- \qnm{Z^l v}^2\right| \leq {\veps}\qnm{Z^lv}^2 + \frac{\lambda}{\veps}\xi(\lambda)\nm{Z^{l+1}v}^2, \quad \forall 0< \veps <1.
\end{equation}
\end{LEM}
\begin{proof}
	1. Since 
	\begin{equation}\label{eq-zv+}
		\nm{Z^l v_+} = \nm{\left(\sum_{k = 0}^s\alpha_kZ^{k}\right)Z^lv } \leq \left(\sum_{k=0}^s |\alpha_k|\lambda^k \right) \nm{Z^l v}= \xi(\lambda)\nm{Z^l v}
	\end{equation} and similarly
	\begin{equation}
		\nm{Z^l v_+ - Z^l v} =  \nm{\left(\sum_{k = 1}^s\alpha_kZ^{k}\right)Z^lv }\leq \left(\sum_{k=1}^s |\alpha_k|\lambda^k \right) \nm{Z^l v} = \lambda \xi(\lambda)\nm{Z^l v},
	\end{equation}
	by using the triangle inequality, one can obtain
	\begin{equation}
	\begin{aligned}
	\left|\nm{Z^l v_+}^2 -  \nm{Z^l v}^2 \right| = &\;\left| \nm{Z^l v_+} - \nm{Z^l v} \right|\left| \nm{Z^l v_+} +  \nm{Z^l v} \right|\\
	\leq &\;  \nm{Z^l v_+-Z^l v} \left| \nm{Z^l v_+} +  \nm{Z^l v} \right|\leq 
	\lambda\xi(\lambda)\nm{Z^l v}^2.
	\end{aligned}
	\end{equation}
		2. We use the definition of $v_+$ to get
	\begin{equation}
	\begin{aligned}
	\qnm{Z^l v_+-Z^l v} =& \qnm{\sum_{k=1}^s\alpha_kZ^{k+l} v}\leq \sum_{k=1}^s|\alpha_k|\qnm{Z^{k-1} \left(Z^{l+1}v\right)}.
	\end{aligned}
	\end{equation}
	We then apply Lemma \ref{lem-Zjump} to obtain 
	\begin{equation}
	\qnm{Z^l v_+-Z^l v}\leq \left(\sqrt{2}\sum_{k=1}^s|\alpha_k|\lambda^{k-\hf}\right)\nm{Z^{l+1} v} = \sqrt{\lambda}\xi(\lambda)\nm{Z^{l+1} v} .
	\end{equation}
	Therefore, with Cauchy--Schwartz inequality and the arithmetic-algebraic inequality $2|ab|\leq \frac{1}{\veps}a^2 + \veps b^2$, one can get
	\begin{equation}
	\begin{aligned}
		\left| \qnm{Z^l v_+}^2- \qnm{Z^l v}^2\right| = &\; \left| \qnm{Z^l v_+- Z^l v + Z^l v}^2- \qnm{Z^l v}^2\right|\\
		= &\; \left| \qnm{Z^l v_+- Z^l v }^2+ 2\qip{Z^l v_+- Z^l v , Z^l v}\right|\\
		\leq &\;	(1+\frac{1}{\veps}) \qnm{Z^l v_+- Z^l v }^2 + \veps \qnm{Z^l v}^2\\
		\leq &\;	 \veps \qnm{Z^l v}^2 + \frac{\lambda}{\veps}\xi(\lambda)\nm{Z^{l+1}v}^2.
	\end{aligned}
	\end{equation}
\end{proof}
Lemma \ref{lem-m-est} will be used in the proof of the modified method.
\begin{LEM}\label{lem-m-est}
	\begin{equation}\label{eq-lemmodi-1}
	\left|\ip{Z^{{l}-1}v,Z^{{l}+1}v} + \nm{Z^{l} v}^2\right|
	\leq
	\veps\qnm{Z^{{l}-1}v}^2+\frac{\lambda}{2\veps}\nm{Z^{l} v}^2.
	\end{equation}
	\begin{equation}\label{eq-lemmodi-2}
	\ip{Z^{{l}-1}v,Z\eta(Z,Z^\top)Z^{{l}+1}v} \leq \frac{\lambda}{\veps} \xi(\lambda)\nm{Z^{l} v}^2 + \veps\qnm{Z^{{l}-1}v}^2,\quad \forall 0<\veps <1.
	\end{equation}
	where $\eta(\cdot,\cdot)$ is a polynomial with two non-communicative variables.
\end{LEM}
\begin{proof} 1. Using discrete integration by parts \eqref{eq-dibp}, one can get 
	\begin{equation}
		\ip{Z^{{l}-1}v,Z^{{l}+1}v} = -\nm{Z^{l} v}^2 - \qip{Z^{{l}-1}v, Z^{l} v}.
	\end{equation}
	 Then we use Cauchy--Schwartz inequality and the  inequality $|ab| \leq \veps a^2 + \frac{1}{4\veps}b^2$ to get
	\begin{equation}\label{eq-f1}
	\begin{aligned}
	\left|\ip{Z^{{l}-1}v,Z^{{l}+1}v} + \nm{Z^{l} v}^2\right|
	\leq \qnm{Z^{{l}-1}v}\qnm{Z^{l} v}
	\leq
	\veps\qnm{Z^{{l}-1}v}^2+\frac{1}{4\veps}\qnm{Z^{l} v}^2.
	\end{aligned}
	\end{equation}
	\eqref{eq-lemmodi-1} can be obtained after substituting the inequality $\qnm{Z^{l} v}\leq \sqrt{2\lambda}\nm{Z^{l}v}$ (obtained from Lemma \ref{lem-Zjump}) into \eqref{eq-f1}.
	
	2. For \eqref{eq-lemmodi-1}, again we first use \eqref{eq-dibp} to get 
	\begin{equation}
	\ip{Z^{{l}-1}v,Z\eta(Z,Z^\top)Z^{{l}+1}v}\\
	=  -\ip{Z^{l} v, \eta(Z,Z^\top)Z^{{l}+1}v}
	-\qip{Z^{{l}-1}v,\eta(Z,Z^\top)Z^{{l}+1}v}.
	\end{equation}
	Then by using a similar estimate as that for \eqref{eq-lemmodi-1}, one can obtain
	\begin{equation}
	\begin{aligned}
	&\ip{Z^{{l}-1}v,Z\eta(Z,Z^\top)Z^{{l}+1}v}+\qip{Z^{{l}-1}v,\eta(Z,Z^\top)Z^{{l}+1}v}\\
	\leq &\;  \nm{\eta(Z,Z^\top)Z}\nm{Z^{l} v}^2 + \qnm{Z^{{l}-1}v}\qnm{\eta(Z,Z^\top)Z^{{l}+1}v}\\
	\leq &\;   \nm{\eta(Z,Z^\top)Z}\nm{Z^{l} v}^2 + 
	\veps\qnm{Z^{k-1}v}^2 + \frac{1}{4\veps}\qnm{\eta(Z,Z^\top)Z^{{l}+1}v}^2\\
	\leq &\; (1+\frac{\lambda}{2\veps})\nm{\eta(Z,Z^\top)Z}\nm{Z^{l} v}^2 + \veps \qnm{Z^{{l}-1}v}^2\\
	\leq &\; \lambda(1+\frac{\lambda}{2\veps})\nm{\eta(Z,Z^\top)}\nm{Z^{l} v}^2 + \veps \qnm{Z^{{l}-1}v}^2. 
	\end{aligned}
	\end{equation}
	The proof is then completed after using the fact $\nm{\eta(Z,Z^\top) }\leq \xi(\lambda)$. 
\end{proof}

\subsection{Proof of the filtering method}
In this section, we prove \eqref{eq-leq} in Lemma \ref{lem-key} for the filtering method. 
\begin{proof}	
	Taking inner-product of \eqref{eq-filter} gives
	\begin{equation}\label{eq-fp}
	\begin{aligned}
	\nm{u_{\rm{F}}}^2 =&\; \nm{u_+}^2 +2\mu\ip{(Z^\top)^{{k^*}-1}Z^{{k^*}}u_+,u_+}+2\nu \ip{(Z^\top)^{{k^*}}Z^{{k^*}}u_+,u_+} \\
	&+ \nm{\left(\mu(Z^\top)^{{k^*}-1}Z^{{k^*}}+\nu(Z^\top)^{{k^*}}Z^{{k^*}}\right)u_+}^2\\
	=&\;\nm{u_+}^2 +2\mu\ip{Z^{{k^*}}u_+,Z^{{k^*}-1}u_+}+2\nu \ip{Z^{{k^*}}u_+,Z^{{k^*}}u_+} \\
	&+ \nm{\left(\mu I +\nu Z^\top\right)(Z^\top)^{{k^*}-1}Z^{{k^*}}u_+}^2\\
	=&\;\nm{u_+}^2 -\mu\qnm{Z^{{k^*}-1}u_+}^2+2\nu \nm{Z^{{k^*}}u_+}^2 
	+ \nm{\left(\mu I+\nu Z^\top\right)(Z^\top)^{{k^*}-1}Z^{{k^*}}u_+}^2.
	\end{aligned}
	\end{equation}
	Here we have used the fact $2\ip{Zv, v} = -\qnm{v}^2$. 
	Invoking Cauchy--Schwartz inequality and  the estimate in \eqref{eq-zv+}, we get
	\begin{equation}
	\begin{aligned}
	\nm{\left(\mu I+\nu Z^\top\right)(Z^\top)^{{k^*}-1}Z^{{k^*}}u_+}
	\leq&\; \nm{\mu I +\nu Z^\top}\nm{Z^\top}^{{k^*}-1}\nm{Z^{{k^*}}u_+}\\
	\leq&\;(\mu + |\nu| \lambda) \lambda^{{k^*}-1}\xi(\lambda)\nm{Z^{{k^*}}u_0}.
	\end{aligned}
	\end{equation}
	Since $(\mu + |\nu| \lambda) \lambda^{{k^*}-1}\xi(\lambda) = \lambda \xi(\lambda)$ under the assumption $k^* = \lceil\frac{p+1}{2}\rceil$ and $\mu = 0$ when $p = 1$, we have
	\begin{equation}\label{eq-fp1}
		\nm{\left(\mu I+\nu Z^\top\right)(Z^\top)^{{k^*}-1}Z^{{k^*}}u_+}^2 \leq \lambda \xi(\lambda) \|Z^{k^*} u_0\|^2. 
	\end{equation}
	On the other hand, we apply Lemma \ref{lem-filter} to obtain 
	\begin{equation}\label{eq-fp2}
	\begin{aligned}
	-\mu\qnm{Z^{{k^*}-1}u_+}^2 \leq &\; -\mu\qnm{Z^{{k^*}-1}u_0}^2 + |\mu|	\left| \qnm{Z^{k^*-1} u_+}^2- \qnm{Z^{k^*-1} u_0}^2\right| \\
	 \leq &\; -\mu\qnm{Z^{{k^*}-1}u_0}^2 + |\mu| \left( {\veps}\qnm{Z^{k^*-1}u_0}^2 + \frac{\lambda}{\veps}\xi(\lambda)\nm{Z^{k^*}u_0}^2\right),
	\end{aligned}
	\end{equation}
	and 
	\begin{equation}\label{eq-fp3}
	\begin{aligned}
	2\nu \nm{Z^{{k^*}}u_+}^2 
	\leq &\; 2\nu \nm{Z^{{k^*}}u_0 }^2 
	+ 2|\nu|	\left|\nm{Z^{k^*} u_+}^2 -  \nm{Z^{k^*} u_0}^2\right|\\
	\leq &\; 2\nu \nm{Z^{{k^*}}u_0 }^2 
	+ |\nu|	\lambda \xi(\lambda)\nm{Z^{k^*} u_0}^2.
	 \end{aligned}
	\end{equation}
	Therefore, one can substitute \eqref{eq-fp1}, \eqref{eq-fp2} and \eqref{eq-fp3} into \eqref{eq-fp} to get
	\begin{equation}
	\begin{aligned}
	\nm{u_{\rm{F}}}^2 \leq &\;\nm{u_+}^2 + (-{\mu}+\veps |\mu| )\qnm{Z^{{k^*}-1}u_0}^2 +\frac{|\mu|}{\veps}\lambda\xi(\lambda)\nm{Z^{{k^*}}u_0}^2 \\
	& + 2\nu \nm{Z^{{k^*}}u_0}^2 + |\nu| \lambda\xi(\lambda)\nm{Z^{{k^*}}u_0}^2\\
	\leq&\; \nm{u_+}^2 + (-{\mu}+\veps|\mu|)\qnm{Z^{{k^*}-1}u_0}^2 + (2\nu +\frac{\lambda}{\veps} \xi(\lambda))\nm{Z^{{k^*}}u_0}^2.
	\end{aligned}
	\end{equation}
	
\end{proof}

\subsection{Proof of the modified method}

Next, we prove Lemma \ref{lem-key} for the modified method. 

\begin{proof}
	Since 
	\begin{equation}
	u_{\rm{M}} = u_+ + \left(u_{\rm{M}} - u_+\right)  = u_+ + \left(R_s(\tZ)-R_s(Z)\right)u_0 = u_+ + \sum_{l=1}^s\alpha_l (\tZ^l - Z^l)u_0,
	\end{equation}
	we then need to evaluate $\tZ^l-Z^l$ for energy estimates. Note that 
	\begin{equation}
	\begin{aligned}
	\tZ^l-Z^l =&\; \left(\left(I+(Z^\top)^{{k^*}-1}(\mu I +\nu Z^\top)Z^{{k^*}-1}\right)^l-I\right)Z^l\\
	=&\; \sum_{k = 1}^l  \binom{l}{k}\left((Z^\top)^{{k^*}-1}(\mu I +\nu Z^\top)Z^{{k^*}-1}\right)^kZ^l\\
	=&\; (Z^\top)^{{k^*}-1}\cdot P_l\cdot Z^{{k^*}+l-1},
	\end{aligned}
	\end{equation}
	where $P_l = (\mu I +\nu Z^\top) Q_l$ and
	\begin{equation}
	\begin{aligned}
	Q_l =&\;l I+Z^{{k^*}-1} \left(\sum_{k = 2}^l  \binom{l}{k}\left((Z^\top)^{{k^*}-1}(\mu I +\nu Z^\top)Z^{{k^*}-1}\right)^{k-2}\right)(Z^\top)^{{k^*}-1}(\mu I +\nu Z^\top).
	\end{aligned}
	\end{equation}
	Therefore, 
	\begin{equation}
	u_{\rm{M}} = u_+ +  (Z^\top)^{{k^*}-1}\sum_{l=1}^s\alpha_lP_lZ^{{k^*}+l-1}u_0.
	\end{equation}
	After taking the inner-product with itself on both sides, one can get
	\begin{equation}\label{eq-m-id}
	\begin{aligned}
	\nm{u_{\rm{M}}}^2 =&\; \nm{u_+}^2 + 2\ip{u_+,\sum_{l=1}^s\alpha_l (Z^\top)^{{k^*}-1}P_lZ^{{k^*}+l-1}u_0} + \nm{(Z^\top)^{{k^*}-1}\sum_{l=1}^s\alpha_l P_lZ^{{k^*}+l-1}u_0}^2\\
	=&\; \nm{u_+}^2 + 2\ip{u_+,\left(\mu Z^{{k^*}-1}+\nu Z^{{k^*}}\right)^\top\sum_{l=1}^s\alpha_l Q_lZ^{{k^*}+l-1}u_0} \\
	&+ \nm{(Z^\top)^{{k^*}-1}\sum_{l=1}^s\alpha_l P_lZ^{{k^*}+l-1}u_0}^2\\
	=&\;\nm{u_+}^2 + 2\ip{\left(\mu Z^{{k^*}-1}+\nu Z^{{k^*}}\right)u_+,\left(\sum_{l=1}^s\alpha_l Q_lZ^{l-1}\right)Z^{{k^*}}u_0} \\
	&+ \nm{(Z^\top)^{{k^*}-1}\left(\sum_{l=1}^s\alpha_l P_lZ^{l-1}\right)Z^{{k^*}}u_0}^2.
	\end{aligned}
	\end{equation}
	We then consider the case $\mu = 0$ and the case $\mu \neq 0$, $p\geq 2$ separately. 
	
	\textbf{Case 1. ($\mu = 0$.)}
	Suppose $\mu = 0$, then using the definition of $u_+$, we have
	\begin{equation}\label{eq-proofm-1}
	\left(\mu Z^{{k^*}-1}+\nu Z^{{k^*}}\right)u_+ = \nu Z^{k^*} u_0 + \xi(Z)Z^{{k^*}+1}u_0.
	\end{equation}
	Since $Q_1 = I$ and $\alpha_1 = 1$, one can get
	\begin{equation}\label{eq-proofm0}
	\left(\sum_{l=1}^s\alpha_l Q_lZ^{l-1}\right)Z^{k^*} u_0 = Z^{k^*} u_0 + \eta(Z,Z^\top)Z^{{k^*}+1} u_0.
	\end{equation}
	Here $\eta(\cdot,\cdot)$ is a polynomial with two non-commuting variables and $\nm{\eta(Z,Z^\top)}\leq \xi(\lambda)$. 
	Therefore, we can expand the inner-product term with \eqref{eq-proofm-1} and \eqref{eq-proofm0}, and invoke Cauchy--Schwartz inequality to obtain
	\begin{equation}\label{eq-m-cross-1}
	2\ip{\left(\mu Z^{{k^*}-1}+\nu Z^{{k^*}}\right)u_+,\left(\sum_{l=1}^s\alpha_l Q_lZ^{l-1}\right)Z^{{k^*}}u_0}\leq (2\nu+\lambda\xi(\lambda)) \nm{Z^{k^*} u_0}^2. 
	\end{equation}
	Furthermore, since $\sum_{l=1}^s\alpha_l P_lZ^{l-1} = \nu Z^\top\left(\sum_{l=1}^s\alpha_l Q_lZ^{l-1}\right)$, we have
	\begin{equation}\label{eq-m-sqres-1}
	\nm{(Z^\top)^{{k^*}-1}\left(\sum_{l=1}^s\alpha_l P_lZ^{l-1}\right)Z^{{k^*}}u_0}^2 =\nm{(Z^\top)^{{k^*}}\eta(Z,Z^\top)Z^{{k^*}}u}^2 \leq \lambda^{2{k^*}}\xi(\lambda)\nm{Z^{{k^*}}u_0}^2.
	\end{equation}
	Since $2k^*\geq 1$, we can substitute \eqref{eq-m-cross-1} and \eqref{eq-m-sqres-1} into \eqref{eq-m-id} to get 
			\begin{equation}
	\nm{u_{\rm{M}}}^2\leq 	\nm{u_+}^2 +(2\nu+\lambda)\nm{Z^{{k^*}}u_0}^2. 
	\end{equation}
	
	\textbf{Case 2.} ($\mu = 1$, $p\geq 2$.)
	The direct evaluation gives
	\begin{equation}\label{eq-proofmodi-1}
	\begin{aligned}
	\left(\mu Z^{{k^*}-1}+\nu Z^{{k^*}}\right)u_+ =&\; \mu Z^{{k^*}-1}u_0 + (\mu+\nu) Z^{{k^*}}u_0+\mu\sum_{l=1}^{s-1}\alpha_{l+1}Z^{{k^*}+l}u_0+\nu\sum_{l=1}^s\alpha_lZ^{{k^*}+l}u_0\\
	=&\; \mu Z^{{k^*}-1}u_0 + (\mu+\nu)Z^{{k^*}}u_0+\xi(Z)Z^{{k^*}+1}u_0.
	\end{aligned}
	\end{equation}
	Since $p\geq 2$, we have $\alpha_2 = \hf$. By using the expression $Q_1 = I$ and $Q_2 = 2I + Z^{k^*-1}(Z^\top)^{k^*-1}(\mu I + \nu Z^\top) = 2I + Z^{k^*-1}\eta(Z,Z^\top)$, one can obtain
	\begin{equation}\label{eq-proofmodi-2}
	\begin{aligned}
	\left(\sum_{l=1}^s\alpha_l Q_lZ^{l-1}\right)Z^{{k^*}}u_0 =&\; Z^{{k^*}}u_0 + Z^{{k^*}+1}u_0+ Z^{{k^*}-1}\eta(Z,Z^\top)Z^{{k^*}+1}u_0.
	\end{aligned}
	\end{equation}
	We can then plug in \eqref{eq-proofmodi-1} and \eqref{eq-proofmodi-2} to expand the inner-product term in \eqref{eq-m-id}. 
	\begin{equation}\label{eq-proofm1}
	\begin{aligned}
	&\;2\ip{(\mu Z^{{k^*}-1}+\nu Z^{{k^*}})u_+,\left(\sum_{l=1}^s\alpha_l Q_lZ^{l-1}\right)Z^{{k^*}}u_0}\\
	=&\;2\mu\ip{Z^{{k^*}-1}u_0,Z^{k^*} u_0}+2(\mu+\nu)\ip{Z^{k^*} u_0,Z^{k^*} u_0} + 2\mu\ip{Z^{{k^*}-1}u_0,Z^{{k^*}+1}u_0} + S \\
	=&\;-\mu\qnm{Z^{k^*} u_0}^2+2\nu\nm{Z^{k^*} u_0}^2 + 2\mu\left(\ip{Z^{{k^*}-1}u_0,Z^{{k^*}+1}u_0} + \nm{Z^{k^*}u_0}^2\right) + S,
	\end{aligned}
	\end{equation}
	where we have used the identity $2\ip{Zv,v} = -\qnm{v}^2$ and we denote by
	\begin{equation}\label{eq-proofm2}
	S := \ip{Z^{{k^*}}u_0,\eta(Z, Z^\top)Z^{{k^*}+1}u_0}+\ip{Z^{{k^*}-1}u_0,Z^{{k^*}-1}\eta(Z,Z^\top)Z^{{k^*}+1}u_0}.
	\end{equation}
	Applying \eqref{eq-lemmodi-1}, one can get
	\begin{equation}\label{eq-proofm3}
	\left|\ip{Z^{{k^*}-1}u_0,Z^{{k^*}+1}u_0} + \nm{Z^{k^*} u_0}^2\right|
	\leq
	\veps\qnm{Z^{{k^*}-1}u_0}^2+\frac{\lambda}{2\veps}\nm{Z^{k^*} u_0}^2,
	\end{equation}
	Note that $p\geq 2$ implies $k^*\geq 2$. Then we use Cauchy--Schwartz inequality to estimate the first term and \eqref{eq-lemmodi-2} for the second term in $S$, which gives 
	\begin{equation}\label{eq-proofm4}
	|S| \leq \frac{\lambda}{\veps} \xi(\lambda)\nm{Z^{k^*} u_0}^2 + \veps\qnm{Z^{{k^*}-1}u_0}^2,\quad \forall 0<\veps <1.
	\end{equation}
   Combining \eqref{eq-proofm1} --\eqref{eq-proofm4}, it can be seen that
	\begin{equation}\label{eq-m-cross-2}
	\begin{aligned}
	&2\ip{(\mu Z^{{k^*}-1}+\nu Z^{{k^*}})u_+,\left(\sum_{l=1}^s\alpha_l Q_lZ^{l-1}\right)Z^{{k^*}}u_0}\\
	\leq&\; (-\mu+\veps|\mu|)\qnm{Z^{{k^*}-1}u_0}^2 + \left(2\nu+ \frac{\lambda}{\veps}\xi(\lambda)\right)\nm{Z^{{k^*}}u_0}^2.
	\end{aligned}
	\end{equation}
	Finally, since $2k^* - 2\geq 1$ for $p \geq 2$, we have
	\begin{equation}\label{eq-m-sqres-2}
	\nm{(Z^\top)^{{k^*}-1}\left(\sum_{l=1}^s\alpha_l P_lZ^{l-1}\right)Z^{{k^*}}u_0}^2\leq \lambda^{2{k^*}-2}\xi(\lambda)\nm{Z^{{k^*}}u_0}^2\leq \lambda \xi(\lambda)\nm{Z^{{k^*}}u_0}^2.
	\end{equation}
	The proof of the modified method can then be completed by substituting \eqref{eq-m-cross-2} and \eqref{eq-m-sqres-2} into \eqref{eq-m-id} . 
\end{proof}

\section{A filtering method for nonlinear problems}\label{sec-nonlinear}
\setcounter{equation}{0}
\setcounter{figure}{0}
\setcounter{table}{0}

For linear problems, Theorem \ref{thm-key} states that the strong stability can be enforced with fixed amount of superviscosity. However, for nonlinear problems \eqref{eq-nlode},
similar analysis no longer goes through. Motivated by the
procedure in \cite{ketcheson2019relaxation}, we pursue the approach
on using the information of the numerical solution to adjust the strength 
of the superviscosity, which leads to the following adaptive filtering method.

\textbf{The filtering method with adaptive $\nu$:} 
Let $u_0$ be the solution at the current time step, and $u_+$ be the solution obtained by the usual RK time integrator. For a given operator $D: V\to V$, if $\nm{D u}>0$, then we apply the filter
\begin{equation}\label{eq-adpfilter}
u_{\rm{F}} = \left(I + \nu D^\top D\right) u_+,\qquad \nu  =  \min\left(\frac{\|u_0\|^2-\|u_+\|^2}{\|D u_+\|^2},0\right).
\end{equation} 
The discussion for linear problems motivates us to choose $D = Z^{k^*}$, while it can be arbitrary in general. 

\begin{PROP} \label{prop-nonlinear}
	Suppose $|\nu |\nm{D}^2 \leq 1$. Then $\nm{u_{\rm{F}}}\leq \nm{u_0}$.
\end{PROP}
\begin{proof}
	If $\nu = 0$, then from \eqref{eq-adpfilter}, we see that $u_{\rm{F}} = u_+$ and $\nm{u_{\rm{F}}} = \nm{u_+}\leq \nm{u_0}$. Otherwise, 	
	$\nu <0$ and $\|u_+\|> \|u_0\|$. By taking inner-product of \eqref{eq-adpfilter} and using Cauchy--Schwartz inequality, it can be shown that,  
	\begin{equation}
	\begin{aligned}
	\nm{u_{\rm{F}}}^2 =&\; \nm{u_+}^2 + 2\nu \nm{D u_+}^2 + \nu^2\nm{D^\top D u_+}^2.
	\end{aligned}
	\end{equation}
	As $|\nu |\nm{D}^2 \leq 1$, we have $|\nu|\nm{D^\top D u_+}^2\leq \nm{D u_+}^2$. Therefore, 
	\begin{equation}
	\begin{aligned}
	\nm{u_{\rm{F}}}^2\leq &\; \nm{u_+}^2 + (2  \nu + |\nu|)\nm{D u_+}^2 \\
	= &\; \nm{u_+}^2 +   \nu \nm{D u_+}^2\\
	=  &\; \nm{u_+}^2 + \frac{\|u_0\|^2-\|u_+\|^2}{\|Du_+\|^2}\nm{D u_+}^2  \\
	= &\; \nm{u_0}^2,
	\end{aligned}
	\end{equation}
	which completes the proof.
\end{proof}
%
%
%
%
%
%

\begin{REM}
	
	Since $u_+ = u(t+\tau) + \mathcal{O}(\tau^{p+1})$ and $\nm{u_+}>\nm{u_0}\geq \nm{u(t + \tau)}$, we have
	\begin{equation}
	\left|\|u_0\|^2 -\|u_+\|^2\right| \leq  \left|\|u(t+\tau)\|^2 -\|u_+\|^2\right| = \mathcal{O}(\tau^{p+1}).
	\end{equation}
	Hence when $\frac{\nm{D^\top D u_+}}{\nm{D u_+}^2} = \mathcal{O}(1)$, the filtering process will not affect the accuracy.
\end{REM}

\begin{REM}
	For linear problems \eqref{eq-ode} with $D = Z^{k^*}$, the post-processing procedure \eqref{eq-adpfilter} retrieves the filtering method \eqref{eq-filter} with $\mu = 0$  and $\nu$ a variant number determined by the solution for each time step. In this case, since $\nm{D u_+}^2 = \mathcal{O}(\tau^{2k^*})$ and the numerator $\nm{u_0}^2 - \nm{u_+}^2 = \mathcal{O}(\tau^{p+1})$, for methods satisfying $k^* = \lceil\frac{p+1}{2}\rceil$, we have $\nu = \mathcal{O}(1)$ if $p$ is odd and $\nu = \mathcal{O}(\tau^{-1})$ if $p$ is even. Note $\nm{D} \leq \lambda^{k^*}$, the assumption in Proposition \ref{prop-nonlinear} can be achieved when $\lambda \ll 1$. (Although we may need $\lambda  = \mathcal{O}(\tau^{\frac{1}{2k^*}})$ if $\nu = \mathcal{O}(\tau^{-1})$, which can be stricter than the usual time step constraint. Similar conditions are also required to achieve weak stability for non-strongly-stable schemes \cite{xu2019L2}. )
\end{REM}
If the mass conservation is not a concern, one can take $D = I$ for stabilization. Since \begin{equation}
\nm{D u_+} =0\Rightarrow \nm{u_+} = 0 \Rightarrow u_{\rm{F}} = u_+ = 0 \Rightarrow\nm{u_{\rm{F}}}\leq \nm{u_0},
\end{equation} 
the strong stability is automatically maintained if $\nm{D^+ u} = 0$.

For ODE systems resulted from semi-discrete schemes of conservation laws, we can take $D$ to be a conservative approximation of $\partial_x^{k^*}$ for preserving the total mass. This retrieves the idea of superviscosity in the context of spectral methods. Note that $\nm{D u_+} = 0$ happens in very rare case. For example, with periodic boundary condition and with Fourier spectral approximation of $\partial_x$, the case $\nm{D u_+}  = 0$ happens only when $u_+$ is a global constant. 

\section{Numerical examples}\label{sec-num}
\setcounter{equation}{0}
\setcounter{figure}{0}
\setcounter{table}{0}

For numerical tests of linear problems, unless otherwise stated, we consider the $p$-stage $p$th order linear RK method 
\begin{equation}
R_p(Z) = \sum_{k=0}^{p}\frac{Z^k}{k!}.
\end{equation}
 The critical values for superviscosity terms are given in Table \ref{tab-critical}.

\begin{table}[h!]
	\centering
	\begin{tabular}{c|c|c|c|c|c|c|c|c|c|c|c|c|c|c|c|c|c|c}
		\hline
		\hline
		$p$&\multicolumn{2}{c|}{$1$}&\multicolumn{2}{c|}{$2$}&\multicolumn{2}{c|}{$3$}&\multicolumn{2}{c|}{$4$}&\multicolumn{2}{c|}{$5$}&\multicolumn{2}{c}{$6$}\\
		\hline
		\hline
		$\mu_0$&\multicolumn{2}{c|}{-}&\multicolumn{2}{c|}{$-\frac{1}{4}$}&\multicolumn{2}{c|}{-}&\multicolumn{2}{c|}{$\frac{1}{144}$}&\multicolumn{2}{c|}{-}&\multicolumn{2}{c}{$-\frac{1}{4800}$}\\
		\hline
		$\nu_0$&\multicolumn{2}{c|}{$-\hf$}&\multicolumn{2}{c|}{$-\frac{1}{8}$}&\multicolumn{2}{c|}{$\frac{1}{24}$}&\multicolumn{2}{c|}{$\frac{1}{144}$}&\multicolumn{2}{c|}{$-\frac{1}{720}$}&\multicolumn{2}{c}{$-\frac{1}{5760}$}\\
		\hline\hline
	\end{tabular}
\caption{Critical values of superviscosity terms for linear RK schemes.}\label{tab-critical}
\end{table}

We also summarize the notations used in Table \ref{tab-notation}.

\begin{table}
	\centering
	\begin{tabular}{c|l}
		\hline
		\hline
		$p$&order of time integrator\\
		\hline
		$k$&order of polynomials in DG discretization\\
		\hline
		$T$&final time\\
		\hline 
		$\tau$&time step size\\
		\hline
		$h$&spatial mesh size\\
		\hline
		$N$&number of spatial mesh cells\\
		\hline
		\multirow{3}{*}{$\alpha$}&parameter for DG discretization; \\
		&$-1$ corresponds to upwind flux;\\
		&$0$ corresponds to central flux  \\
		\hline
		$(\mu,\nu)$&coefficients for superviscosity\\
		\hline
		$u_h$&fully discrete solution with $(\mu, \nu) = 0$\\
		\hline
		$u_{h,\rm{M}}$&fully discrete solution with the modified method\\
		\hline
		$u_{h,\rm{F}}$&fully discrete solution with the filtering method\\
		\hline
		\hline
	\end{tabular}
	\caption{Chart of notations}\label{tab-notation}
\end{table}
\subsection{A non-normal linear ODE test}\label{examp-ode}
In this numerical test, we consider a linear autonomous ODE system
$
\frac{d}{dt}u = Lu,
$ with non-normal right hand side introduced in \cite{sun2017rk4}, 
\begin{equation}
	L =- \left(\begin{matrix}
	1&2&2\\
	0&1&2\\
	0&0&1
	\end{matrix}\right).
\end{equation}

\begin{table}
	\small
	\begin{tabular}{c|c|c|c|c|c|c|c|c|c}
		\hline
		\hline
		$p$&$\mu$&$\nu$&$\tau$&$10^{-1}$&$10^{-2}$&$10^{-3}$&$10^{-4}$&$10^{-5}$&$10^{-6}$\\
		\hline
		\hline
		\multirow{4}[7]{*}{$1$}&$0$&$0$&&\ 1.72E-02&\ 1.52E-04&\ 1.50E-06&\ 1.50E-08&\ 1.50E-10&\ 1.50E-12\\
		\cline{2-10}
		&\multirow{2}{*}{$0$}&\multirow{2}{*}{$-\hf$}
		&M&\ 4.64E-03&\ 4.04E-06&\ 4.00E-09&\ 4.00E-12&\ 4.00E-15&\ 4.00E-18\\
		&&&F&\ 7.09E-03&\ 7.88E-06&\ 7.99E-09&\ 8.00E-12&\ 8.00E-15&\ 8.00E-18\\
		\cline{2-10}
		&\multirow{2}{*}{$0$}&\multirow{2}{*}{$-\frac{1.01}{2} $}
		&M&\ 4.53E-03&\ 2.58E-06&-1.67E-09&-1.67E-11&-1.67E-13&-1.67E-15\\
		&&&F&\ 7.00E-03&\ 6.46E-06&-1.67E-09&-1.67E-11&-1.67E-13&-1.67E-15\\
		\hline
		\hline
		\multirow{4}[16]{*}{$2$}&$0$&$0$&&\ 1.44E-05&\ 1.27E-09&\ 1.25E-13&\ 1.25E-17&\ 1.25E-21&\ 1.25E-25\\
		\cline{2-10}
		&\multirow{2}{*}{$-\frac{1}{4}$}&\multirow{2}{*}{$-\frac{1}{8}$}&M&\ 6.37E-06&\ 5.64E-10&\ 5.52E-14&\ 5.51E-18&\ 5.51E-22&\ 5.51E-26\\
		&&&F&\ 1.27E-05&\ 1.12E-09&\ 1.10E-13&\ 1.10E-17&\ 1.10E-21&\ 1.10E-25\\
		\cline{2-10}
		&\multirow{2}{*}{$-\frac{0.99}{4}$}&\multirow{2}{*}{$-\frac{1.01}{8}$}&M&\ 6.17E-06&\ 4.24E-10&\ 1.21E-14&\ 4.33E-20&-1.08E-23&-1.23E-27\\
		&&&F&\ 1.24E-05&\ 9.57E-10&\ 4.08E-14&\ 5.10E-19&-5.81E-24&-1.18E-27\\
		\cline{2-10}
		&\multirow{2}{*}{$0$}&\multirow{2}{*}{$-\frac{1.01}{8}$}&M&\ 2.80E-06&\ 2.05E-11&-9.14E-16&-1.22E-19&-1.25E-23&-1.25E-27\\
		&&&F&\ 5.22E-06&\ 5.29E-11&-5.79E-16&-1.18E-19&-1.24E-23&-1.25E-27\\
		\hline
		\hline
		\multirow{4}[7]{*}{$3$}&$0$&$0$&&-5.14E-06&-4.26E-10&\ 4.18E-14&\ 4.17E-18&\ 4.17E-22&-4.17E-26\\
		\cline{2-10}
		&\multirow{2}{*}{$0$}&\multirow{2}{*}{$\frac{1}{24}$}&M&-1.04E-07&-3.54E-12&-3.69E-17&-3.70E-22&-2.70E-27&-3,70E-32\\
		&&&F&-1.54E-06&-1.49E-11&-1.48E-16&-1.48E-21&-1.48E-26&-1.48E-31\\
		\cline{2-10}
		&\multirow{2}{*}{$0$}&\multirow{2}{*}{$\frac{1.01}{24}$}&M&-4.24E-08&\ 7.40E-13&\ 3.81E-16&\ 4.13E-20&\ 4.16E-24&\ 4.17E-28\\
		&&&F&-1.50E-06&-1.07E-11&\ 2.68E-16&\ 4.02E-20&\ 4.15E-24&\ 4.17E-28\\
		\hline
		\hline
		\multirow{4}[16]{*}{$4$}&$0$&$0$&&\ 2.22E-07&\ 1.89E-12&\ 1.86E-17&\ 1.85E-22&\ 1.85E-27&\ 1.85E-32\\
		\cline{2-10}
		&\multirow{2}{*}{$\frac{1}{144}$}&\multirow{2}{*}{$\frac{1}{144}$}&M&\ 1.13E-07&\ 1.11E-13&\ 1.11E-19&\ 1.11E-25&\ 1.11E-31&\ 1.11E-37\\
		&&&F&\ 1.48E-07&\ 1.65E-13&\ 1.66E-19&\ 1.67E-25&\ 1.67E-31&\ 1.67E-37\\
		\cline{2-10}
		&\multirow{2}{*}{$\frac{1.01}{144}$}&\multirow{2}{*}{$\frac{0.99}{144}$}&M&\ 1.10E-07&\ 9.20E-14&-7.48E-20&-1.74E-24&-1.84E-29&-1.85E-34\\
		&&&F&\ 1.46E-07&\ 1.46E-13&-1.89E-20&-1.69E-24&-1.84E-29&-1.85E-34\\
		\cline{2-10}
		&\multirow{2}{*}{$\frac{1.01}{144}$}&\multirow{2}{*}{$0$}&M&\ 3.48E-08&\ 1.64E-14&-1.50E-19&-1.82E-24&-1.85E-29&-1.85E-34\\
		&&&F&\ 8.59E-08&\ 7.20E-14&-9.44E-20&-1.76E-24&-1.84E-29&-1.85E-34\\
		\cline{2-10}
		&\multirow{2}{*}{$0$}&\multirow{2}{*}{$-100$}&M&-7.67E-05&-1.07E-09&-1.08E-15&-9.15E-22&\ 7.52E-28&\ 1.74E-32\\
		&&&F&-7.18E-05&-1.05E-09&-1.08E-15&-9.15E-22&\ 7.52E-28&\ 1.74E-32\\
		\hline
		\hline
	\end{tabular}
\caption{$\|R_p(Z)\|-1$ in Section \ref{examp-ode2} with different magnitude of superviscosity $(\mu,\nu)$ and different time step sizes $\tau$. $\rm{M}$ corresponds to the modified method and $\rm{F}$ corresponds to the filtering method. }\label{tab-ode-norm}
\end{table}

\begin{table}[h!]
	\small
	\centering
	\begin{tabular}{ c | c | c | c | c |c | c }
		\hline
		\hline
		\multirow{2}{*}{$p$}&\multirow{2}{*}{$(\mu,\nu)$}&\multirow{2}{*}{$\tau$}&\multicolumn{2}{c|}{$u_{\rm{M}}$}&\multicolumn{2}{c}{$u_{\rm{F}}$}\\
		\cline{4-7}
		&&&$l^2$ error &order&$l^2$ error &order\\
		\hline
		\hline
		\multirow{5}{*}{$1$}&
		\multirow{5}{*}{$(0,-1)$}
		&    1/20&1.9177E-01&-   &1.6907E-01&-\\
		&&   1/40&1.0328E-01&0.89&9.6610E-02&0.81\\
		&&   1/80&5.3686E-02&0.94&5.1868E-02&0.90\\
		&&  1/160&2.7384E-02&0.97&2.6908E-02&0.95\\
		&&  1/320&1.3831E-02&0.99&1.3709E-02&0.97\\
		\hline\hline
		\multirow{5}{*}{$2$}&
		\multirow{5}{*}{$(1,-1)$}
		&    1/20&3.7343E-02&-	 &3.5741E-02&-\\
		&&   1/40&8.9226E-03&2.07&8.7225E-03&2.03\\
		&&   1/80&2.1629E-03&2.04&2.1383E-03&2.02\\
		&&  1/160&5.3140E-04&2.03&5.2838E-04&2.02\\
		&&  1/320&1.3641E-04&2.01&1.3126E-04&2.01\\
		\hline\hline
		\multirow{5}{*}{$3$}&
		\multirow{5}{*}{$(0,-1)$}
		&1/20&5.2177E-03&-   &4.9703E-03&-\\
		&&   1/40&6.5467E-04&2.99&6.3896E-04&2.96\\
		&&   1/80&8.1878E-05&3.00&8.0892E-05&2.98\\
		&&  1/160&1.0236E-05&3.00&1.0174E-05&2.99\\
		&&  1/320&1.2795E-06&3.00&1.2756E-06&3.00\\
		\hline\hline
		\multirow{5}{*}{$3$}&
		\multirow{5}{*}{$(1,0)$}
		&	 1/20&3.1195E-02&-	 &2.9831E-02&-\\
		&&   1/40&7.9716E-03&1.97&7.7898E-03&1.94\\
		&&   1/80&2.0042E-03&1.99&1.9810E-03&1.98\\
		&&  1/160&5.0178E-04&2.00&4.9888E-04&1.99\\
		&&  1/320&1.2550E-04&2.00&1.2513E-04&2.00\\
		\hline\hline
		\multirow{5}{*}{$4$}&
		\multirow{5}{*}{$(1,-1)$}
		&    1/20&5.9784E-04&-	 &5.7407E-04&-\\	
		&&   1/40&3.5563E-05&4.07&3.4852E-05&4.04\\
		&&   1/80&2.1659E-06&4.04&2.1442E-06&4.02\\
		&&  1/160&1.3360E-07&4.02&1.3293E-07&4.01\\
		&&  1/320&8.2943E-09&4.01&8.2735E-09&4.01\\
		\hline
		\hline
	\end{tabular}
\caption{Accuracy test in Section \ref{examp-ode-acc} for the modified and filtering methods with a non-normal ODE system.}\label{tab-ode-accuracy}
\end{table}

\subsubsection{Norm of RK operators}\label{examp-ode2}
We first examine the norm of RK operators with different magnitude of superviscosity. As we can see from Table \ref{tab-ode-norm}, the original first, second and fourth  order 
RK methods can not ensure $\|R_p(Z)\|\leq 1$ under small time steps and are hence not strongly stable. While after adding
the superviscosity terms with $\mu>\mu_0$ and $\nu<\nu_0$, even though the deviation of $\mu$ and $\nu$ is only one percent away from the critical values, we have $\|R_p(Z)\|\leq 1$ with $\tau$ sufficiently small. Both the modified and the filtering approach can stabilize the methods and their difference is quite small, especially when $\tau$ approaches zero. It should also be noted 
that, for the fourth order methods, if $\mu = 0$, even when $\nu$ is as negative as $-100$, the method would not be strongly stable for small $\tau$.  This validates the necessity to include dispersive superviscosity. The original third order RK method itself is strongly stable. If we add an anti-superviscosity term to the method, the method is still strongly stable if $\nu\leq \nu_0$, which also indicates critical values given in Section \ref{sec-mf} are sharp.

\subsubsection{Accuracy test}\label{examp-ode-acc}
Then we test the accuracy of the modified and the filtering RK methods. The initial condition is taken as $u(0) = (1\  1\  1)^T$ and the exact solution of the initial value problem is 
\begin{equation}
u = \left(\begin{array}{c}
1-4t+2t^2\\
1-2t\\
1
\end{array}\right)e^{-t}.
\end{equation}
In Table \ref{tab-ode-accuracy}, we document the $l^2$ error of the solution. It can be seen that, when $\nu = \mathcal{O}(1)$, the modified and the filtering RK methods maintain the original order of accuracy. If we have $\mu = \mathcal{O}(1)$, methods of even order still achieve the designed convergence rate. But the odd order methods may have a degenerated order of accuracy. This confirms our requirement $\mu = 0$ for odd $p$.

\subsection{DG method for linear advection equation}
In this test, we study the modified and the filtering methods for time integration of the linear advection equation
\begin{equation}\label{eq-adv}
	u_t + u_x = 0,\qquad t>0, \qquad x\in (0,2\pi)
\end{equation}
with the periodic boundary condition. DG method is used for spatial discretization. Here we only introduce the most basic concepts that are needed to explain our numerical scheme. For the systematic explanation of the method, we refer to a series of paper by Cockburn et al.  \cite{rkdg1,rkdg2,rkdg3,rkdg4,rkdg5}. 

We consider a uniform partition of the spatial domain $\Omega = (0,2\pi) = \cup_{j = 1}^NI_j$ with the mesh size $h$. Let $P^k(I_j)$ be the linear space spanned by polynomials on $I_j$ with degree no larger than $k$. We take the discontinuous polynomial space 
\begin{equation}
V_h = \{v\in L^2(\Omega): v\big|_{I_j} \in P^k(I_j)\}
\end{equation}
as the finite element space. Since functions in $V_h$ can be double-valued at the cell interfaces, we use $v^+$ and $v^-$ to represent the right and left limits of $v$.

For the DG approximation to \eqref{eq-adv}, we look for $u_h \in V_h$, such that for any $v\in V_h$,
\begin{equation}
\int_{I_j} (u_h)_t v_h dx  =  \int_{I_j} u_h v_x dx - \widehat{u}_{j+\hf}v^-_{j+\hf} + \widehat{u}_{j-\hf}v^+_{j-\hf}.
\end{equation}
Here $\widehat{u} =  \frac{1-\alpha}{2}u^-_h + \frac{1+\alpha}{2} u_h^+$ is the numerical flux. $\alpha = -1$ corresponds to the upwind flux and $\alpha = 0$ corresponds to the central flux. 

We introduce the operator $L_\alpha : V_h \to V_h$, that is uniquely determined by the variational form 
\begin{equation}
\int_{I_j} L_\alpha u_h v_h dx  =  \int_{I_j} u_h v_x dx - \widehat{u}_{j+\hf}v^-_{j+\hf} + \widehat{u}_{j-\hf}v^+_{j-\hf}, \qquad \forall v_h \in V_h,
\end{equation}
and we consider the $L^2$ inner-product on $V_h$. Then it can be shown that
\begin{equation}
\ip{L_\alpha v, v} = \frac{\alpha}{2}\sum_{j = 1}^N \left(v_{j+\hf}^+-v_{j+\hf}^-\right)^2\leq 0, \qquad \text{ if } \alpha \leq 0,
\end{equation}
and (see, for example \cite{sun2017stability} and \cite{zhang2012fully})
\begin{equation}
	L_\alpha^\top  = -L_{-\alpha}.
\end{equation}
For $\alpha \leq 0$, the semi-discrete scheme satisfies
\begin{equation}
\frac{d}{dt}u_h = L_\alpha u_h, \qquad \ip{L_\alpha v, v}  \leq 0, \qquad \forall v\in V_h.
\end{equation}

\subsubsection{Norm of the RK operators}\label{examp-adv}
Again, we begin with testing the norm of RK operators. It can be shown that, the operator $Z$ does not explicitly depend on the spatial domain $\Omega$. Instead, it relates with the order $p$, polynomial degree $k$, the number of spatial mesh $N$ and the CFL number $\frac{\tau}{h}$. We fix $N = 10$ in the test. The detailed examination is given in Table \ref{tab-adv-norm}. One can see from the table, with the prescribed ranges of superviscosity in Section \ref{sec-mf}, both the modified and the filtering methods are strongly stable for sufficiently small $\tau$.

\begin{table}
	\small
	\begin{tabular}{c|c|c|c|c|c|c|c|c|c|c}
		\hline
		\hline
		$p$&$k$&$\mu$&$\nu$&${\tau}/{h}$&$10^{-1}$&$10^{-2}$&$10^{-3}$&$10^{-4}$&$10^{-5}$&$10^{-6}$\\
		\hline
		\hline
		\multirow{4}[7]{*}{$1$}&	\multirow{4}[7]{*}{$1$}&$0$&$0$&& 6.44E-02& 6.06E-04& 6.01E-06& 6.00E-08& 6.00E-10& 6.00E-12\\
		\cline{3-11}
		&&\multirow{2}{*}{$0$}&\multirow{2}{*}{$-\frac{1}{2} $}&M& 1.47E-02& 1.22E-05& 1.20E-08& 1.20E-11&1.20E-14& 1.20E-17\\
		\cline{5-11}
		&&&&F& 1.64E-02& 2.33E-05& 2.39E-08& 2.40E-11& 2.40E-14& 2.40E-17\\
		\cline{3-11}
		&&\multirow{2}{*}{$0$}&\multirow{2}{*}{$-\frac{1.01}{2} $}&M& 1.43E-02& 6.37E-06& 0.00E-99& 0.00E-99& 0.00E-99& 0.00E-99\\
		\cline{5-11}
		&&&&F& 1.59E-02& 1.75E-05& 0.00E-99& 0.00E-99& 0.00E-99& 0.00E-99\\
		\hline
		\hline
		\multirow{4}[7]{*}{$2$}&	\multirow{4}[7]{*}{$2$}&$0$&$0$&& 1.85E-03& 1.52E-07& 1.50E-11& 1.50E-15& 1.50E-19& 1.50E-23\\
		\cline {3-11}
		&&\multirow{2}{*}{$-\frac{1}{4}$}&\multirow{2}{*}{$-\frac{1}{8} $}&M& 0.00E-99& 1.29E-08& 1.24E-12& 1.23E-16& 1.23E-20& 1.23E-24\\
		\cline{5-11}
		&&&&F& 4.08E-04& 2.62E-08& 2.47E-12& 2.46E-16& 2.46E-20& 2.46E-24\\
		\cline {3-11}
		&&\multirow{2}{*}{$-\frac{0.99}{4}$}&\multirow{2}{*}{$-\frac{1.01}{8} $}&M& 0.00E-99& 1.08E-08& 5.56E-13& 0.00E-99& 0.00E-99& 0.00E-99\\
		\cline{5-11}
		&&&&F& 3.92E-04& 2.39E-08& 1.63E-12& 3.39E-17& 0.00E-99& 0.00E-99\\
		\cline{3-11}
		&&\multirow{2}{*}{$0$}&\multirow{2}{*}{$-\frac{1.01}{8} $}&M& 6.09E-04& 1.30E-09& 0.00E-99& 0.00E-99& 0.00E-99& 0.00E-99\\
		\cline{5-11}
		&&&&F& 3.62E-04& 3.94E-09& 0.00E-99& 0.00E-99& 0.00E-99& 0.00E-99\\
		\hline
		\hline
			\multirow{4}[7]{*}{$3$}&	\multirow{4}[7]{*}{$3$}&$0$&$0$&& 3.65E-02& 0.00E-99& 0.00E-99& 0.00E-99& 0.00E-99& 0.00E-99\\
		\cline{3-11}
		&&\multirow{2}{*}{$0$}&\multirow{2}{*}{$\frac{1}{24} $}&M& 6.40E-02& 1.81E-11& 1.66E-17& 1.65E-23& 1.65E-29& 1.65E-35\\
		\cline{5-11}
			&&&&F& 1.09E-01& 1.72E-11& 1.65E-17& 1.65E-23& 1.65E-29& 1.65E-35\\
		\cline {3-11}
		&&\multirow{2}{*}{$0$}&\multirow{2}{*}{$\frac{0.99}{24} $}&M& 6.29E-02& 0.00E-99& 0.00E-99& 0.00E-99& 0.00E-99& 0.00E-99\\
		\cline{5-11}
		&&&&F& 1.07E-01& 0.00E-99& 0.00E-99& 0.00E-99& 0.00E-99& 0.00E-99\\
		\hline
		\hline
		\multirow{4}[7]{*}{$4$}&	\multirow{4}[7]{*}{$4$}&$0$&$0$&& 2.82E-01& 1.70E-07& 1.47E-12& 1.45E-17& 1.45E-22& 1.45E-27\\
		\cline{3-11}
		&&\multirow{2}{*}{$\frac{1}{144}$}&\multirow{2}{*}{$\frac{1}{144} $}&M& 7.94E-01& 8.93E-08& 8.71E-14& 8.68E-20& 8.68E-26& 8.68E-32\\
		\cline{5-11}
		&&&&F& 2.78E+00& 1.17E-07& 1.29E-13& 1.30E-19& 1.30E-25& 1.30E-31\\
		\cline{3-11}
		&&\multirow{2}{*}{$\frac{1.01}{144}$}&\multirow{2}{*}{$\frac{0.99}{144} $}&M& 7.38E-01&8.72E-08& 7.20E-14& 6.03E-22& 0.00E-99& 0.00E-99\\
		\cline{5-11}
		&&&&F& 2.72E+00& 1.15E-07& 1.14E-13& 2.13E-22& 2.92E-30& 0.00E-99\\
		\cline {3-11}
		&&\multirow{2}{*}{$\frac{1.01}{144}$}&\multirow{2}{*}{$0 $}&M& 6.91E+00& 2.43E-08& 8.68E-15& 0.00E-99& 0.00E-99& 0.00E-99\\
		\cline{5-11}
		&&&&F& 2.05E+01& 6.38E-08& 5.21E-14& 0.00E-99& 0.00E-99& 0.00E-99\\
		\hline
		\hline
		\multirow{4}[7]{*}{$5$}&	\multirow{4}[7]{*}{$5$}&$0$&$0$&& 2.20E+00& 2.45E-09& 2.26E-15& 2.24E-21& 2.23E-27& 2.23E-33\\
		\cline {3-11}
		&&\multirow{2}{*}{$0$}&\multirow{2}{*}{$-\frac{1}{720} $}&M& 5.14E+02& 3.83E-10& 4.79E-17& 4.94E-24& 4.95E-31& 4.96E-38\\
		\cline{5-11}
		&&&&F& 1.20E+01& 5.46E-10& 7.29E-17& 7.60E-24& 7.63E-31& 7.63E-38\\
		\cline{3-11}
		&&\multirow{2}{*}{$0$}&\multirow{2}{*}{$-\frac{1.01}{720}$}&M& 5.31E+02& 3.63E-10& 2.60E-17& 0.00E-93& 0.00E-92& 0.00E-91\\
		\cline{5-11}
		&&&&F& 1.22E+01& 5.27E-10& 5.13E-17& 0.00E-99& 0.00E-99& 0.00E-99\\
				\hline
		\hline
		\multirow{4}[7]{*}{$6$}&	\multirow{4}[7]{*}{$6$}&$0$&$0$&& 9.76E+00& 2.73E-12& 2.53E-20& 2.52E-28& 2.52E-36& 2.52E-44\\
		\cline {3-11}
		&&\multirow{2}{*}{$-\frac{1}{4800}$}&\multirow{2}{*}{$-\frac{1}{5760} $}&M& 2.30E+08& 4.07E-13& 2.79E-21& 2.69E-29& 2.67E-37& 2.67E-45\\
		\cline{5-11}
		&&&&F& 5.74E+02& 6.54E-13& 4.47E-21& 4.30E-29& 4.28E-37& 4.28E-45\\
		\cline{3-11}
		&&\multirow{2}{*}{$-\frac{0.99}{4800}$}&\multirow{2}{*}{$-\frac{1.01}{5760}$}&M& 2.53E+08& 3.83E-13& 2.33E-21& 1.08E-29& 0.00E-99& 0.00E-99\\
		\cline{5-11}
		&&&&F& 5.84E+02& 6.30E-13& 3.99E-21& 2.36E-29& 3.16E-38& 0.00E-99\\
				\cline{3-11}
		&&\multirow{2}{*}{$0$}&\multirow{2}{*}{$-\frac{1.01}{5760}$}&M& 1.26E+09& 2.82E-13& 1.52E-22& 0.00E-99& 0.00E-99& 0.00E-99\\
		\cline{5-11}
		&&&&F& 7.77E+02& 4.61E-13& 4.26E-22& 0.00E-99& 0.00E-99& 0.00E-99\\
		\hline
		\hline
	\end{tabular}
	\caption{$\|R_p(Z)\|-1$ in Section \ref{examp-adv} with different magnitude of superviscosity $(\mu,\nu)$ and different CFL number $\tau/h$. The number of spatial cell is set as $N = 10$. $\rm{M}$ corresponds to the modified method and $\rm{F}$ corresponds to the filtering method. }\label{tab-adv-norm}
\end{table}

\subsubsection{Accuracy test}\label{examp-adv-acc}

Then we examine the accuracy of the fully discrete scheme with modified and filtering time integrators. We solve the linear advection equation \eqref{eq-adv} with the initial condition $u(x,0) = e^{\sin(x)}$. The upwind DG method ($\alpha = -1$) with polynomials of degree $k = p-1$ is used for the spatial discretization. We compute with the time step $\tau = 0.02 h$ to the final time $T = 1$. The difference between the modified and the filtering method is quite negligible. They both achieves the designed convergence rate. 

 \begin{table}[h!]
 	\small
	\centering
	\begin{tabular}{ c | c | c | c | c |c | c |c}
		\hline\hline
		\multirow{2}{*}{$p$}&\multirow{2}{*}{$k$}&\multirow{2}{*}{$(\mu,\nu)$}&\multirow{2}{*}{$N$}&\multicolumn{2}{c|}{$u_{\rm{M}}$} &\multicolumn{2}{c}{$u_{\rm{F}}$} \\ 
		\cline{5-8}
		 &&&&$L^2$ error &order&$L^2$ error &order\\
\hline\hline
\multirow{5}{*}{$1$}&
\multirow{5}{*}{$0$}&
\multirow{5}{*}{$(0,-1)$}
&   20&4.2660E-01&-&4.2924E-01& - \\
&&&   40&2.3032E-01& 0.89&2.3196E-01& 0.89 \\
&&&   80&1.2005E-01& 0.94&1.2096E-01& 0.94 \\
&&&  160&6.1391E-02& 0.97&6.1574E-02& 0.97 \\
&&&  320&3.1042E-02& 0.98&3.1058E-02& 0.99 \\
 \hline\hline
 \multirow{5}{*}{$2$}&
\multirow{5}{*}{$1$}&
\multirow{5}{*}{$(1,-1)$}
 &   20&1.7236E-02& -&1.7224E-02& - \\
&&&   40&4.3120E-03& 2.00&4.3109E-03& 2.00 \\
&&&   80&1.0761E-03& 2.00&1.0760E-03& 2.00 \\
&&&  160&2.8280E-04& 1.93&2.8280E-04& 1.93 \\
&&&  320&7.0714E-05& 2.00&7.0713E-05& 2.00 \\
 \hline\hline
\multirow{5}{*}{$3$}&
\multirow{5}{*}{$2$}&
\multirow{5}{*}{$(0,-1)$}
 &   20&7.3661E-04& -&7.3662E-04& - \\
&&&   40&9.1956E-05& 3.00&9.1957E-05& 3.00 \\
&&&   80&1.1436E-05& 3.01&1.1437E-05& 3.01 \\
&&&  160&1.5581E-06& 2.88&1.5581E-06& 2.88 \\
&&&  320&1.9471E-07& 3.00&1.9471E-07& 3.00 \\
 \hline\hline
 \multirow{5}{*}{$4$}&
\multirow{5}{*}{$3$}&
\multirow{5}{*}{$(1,-1)$}
 &   20&3.0535E-05& -&3.0538E-05& -\\
&&&   40&1.9265E-06& 3.99&1.9265E-06& 3.99 \\
&&&   80&1.1785E-07& 4.03&1.1785E-07& 4.03 \\
&&&  160&8.4115E-09& 3.81&8.4115E-09& 3.81 \\
&&&  320&5.2557E-10& 4.00&5.2557E-10& 4.00 \\
 \hline\hline
\multirow{5}{*}{$5$}&
\multirow{5}{*}{$4$}&
\multirow{5}{*}{$(0,-1)$}
 &   20&1.1948E-06& - &1.1949E-06& - \\
&&&   40&3.8144E-08& 4.97 &3.8143E-08& 4.97 \\
&&&   80&1.1551E-09& 5.05&1.1551E-09& 5.05 \\
&&&  160&4.2732E-11& 4.76&4.2732E-11& 4.76 \\
&&&  320&1.3305E-12& 5.01&1.3305E-12& 5.01 \\
\hline\hline
	\end{tabular}
\caption{Accuracy test in Section \ref{examp-adv-acc} with upwind DG discretization for \eqref{eq-adv} with modified and filtering methods.}
\end{table}

\subsubsection{Energy dissipation and compensation}

We now use particular initial value problems to illustrate how the extra superviscosity terms would effect the monotonicity of the $L^2$ norm of the solution. 

\begin{figure}[h!]
	\centering
	\begin{subfigure}{0.32\textwidth}
		\includegraphics[width=\textwidth]{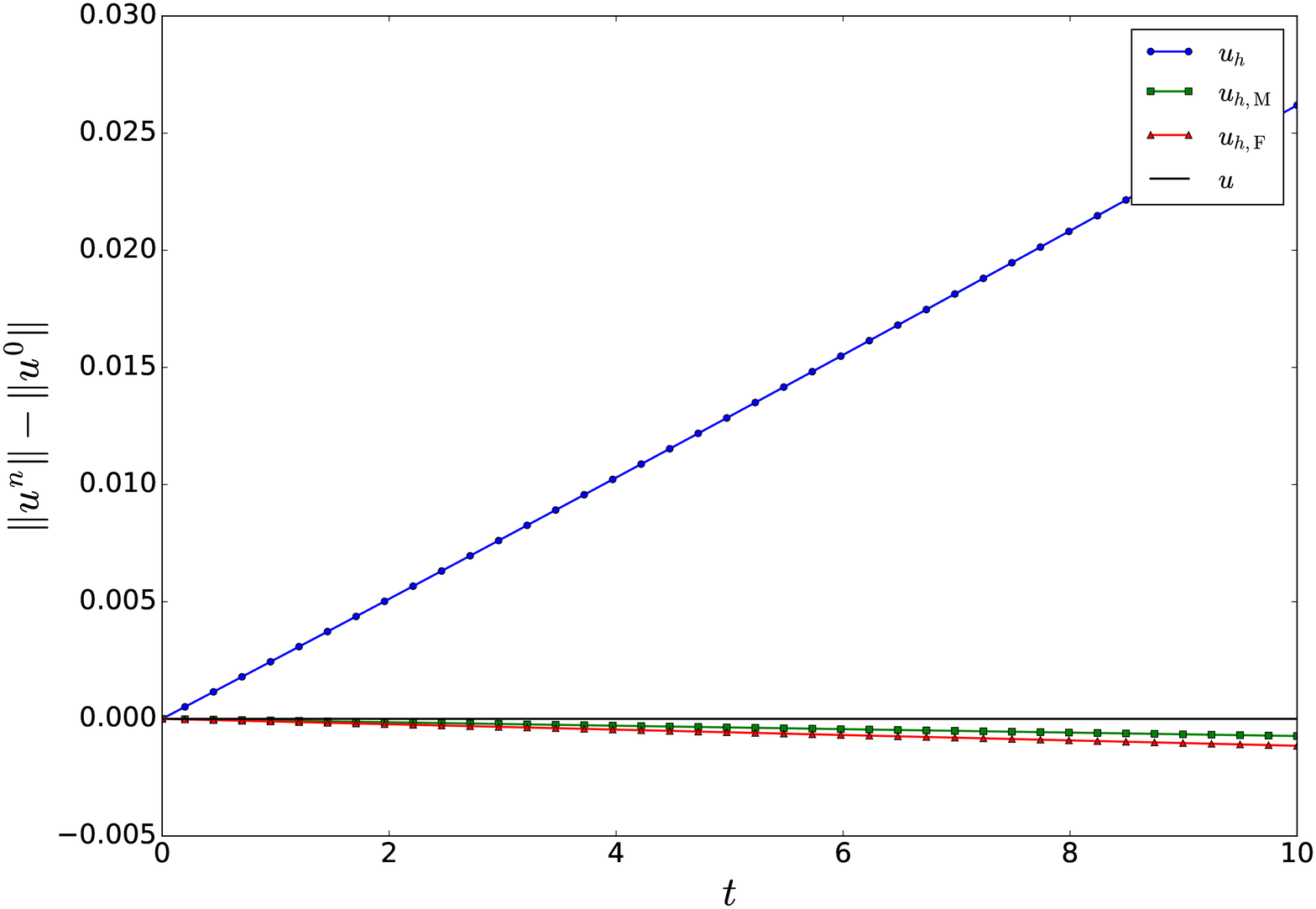}
		\caption{$p = k = 1$, $\nu = -\frac{1.01}{2}$.}	
	\end{subfigure}
\begin{subfigure}{0.32\textwidth}
	\includegraphics[width=\textwidth]{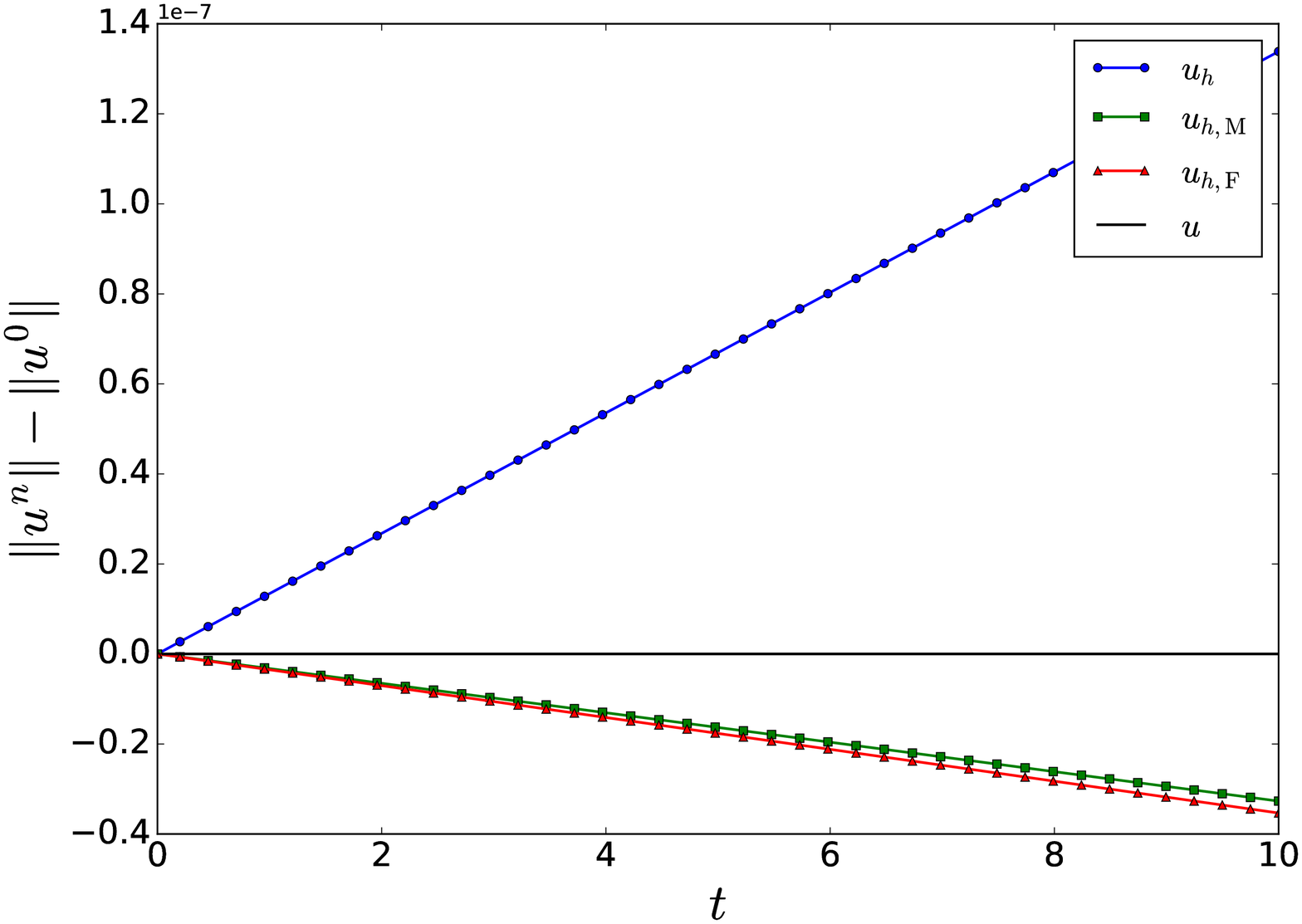}
	\caption{$p = k = 2$, $\nu = -\frac{1.01}{8}$.}	
\end{subfigure}
	\begin{subfigure}{0.32\textwidth}
		\includegraphics[width=\textwidth]{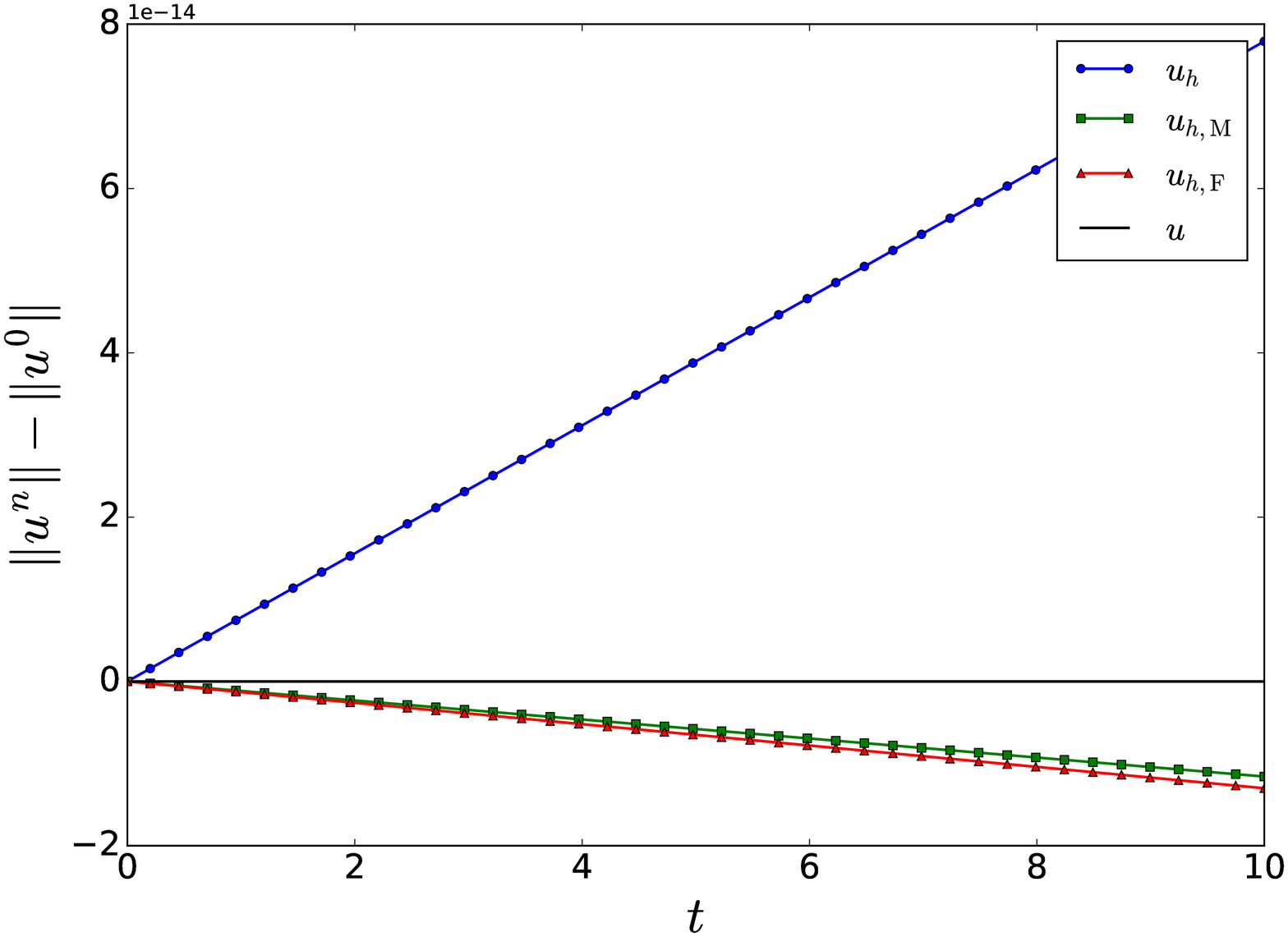}
		\caption{$p = k = 5$,
			 $\nu = -\frac{1.1}{720}$.}	
	\end{subfigure}
\caption{Relative $L^2$ norm of the solutions with superviscosity. $\mu = 0$. }\label{fig-energy-decay}
\end{figure}

For the $p$th order methods with $p = 1,2, 5$, when coupled with $P^p$ 
upwind DG spatial discretization, the solution may have an increasing $L^2$ norm and may not even be stable for $p = 1,2$ under the usual CFL condition. We attempt to enforce the strong stability by including a superviscosity term. We set the initial condition to be $u(x,0) = e^{\sin{x}}$ and use $80$-cell upwind DG method for the
spatial discretization. The time step is set as $\tau = 0.05h$ and we compute to $T = 10$. As we can see from Figure \ref{fig-energy-decay}, the original time discretization methods produce solutions with an increasing $L^2$ norm. While by adding a small superviscosity as that indicated in Section \ref{sec-mf}, the $L^2$ norm becomes decaying. 

\begin{figure}[h!]
	\centering
	\begin{subfigure}{0.45\textwidth}
		\includegraphics[width=\textwidth]{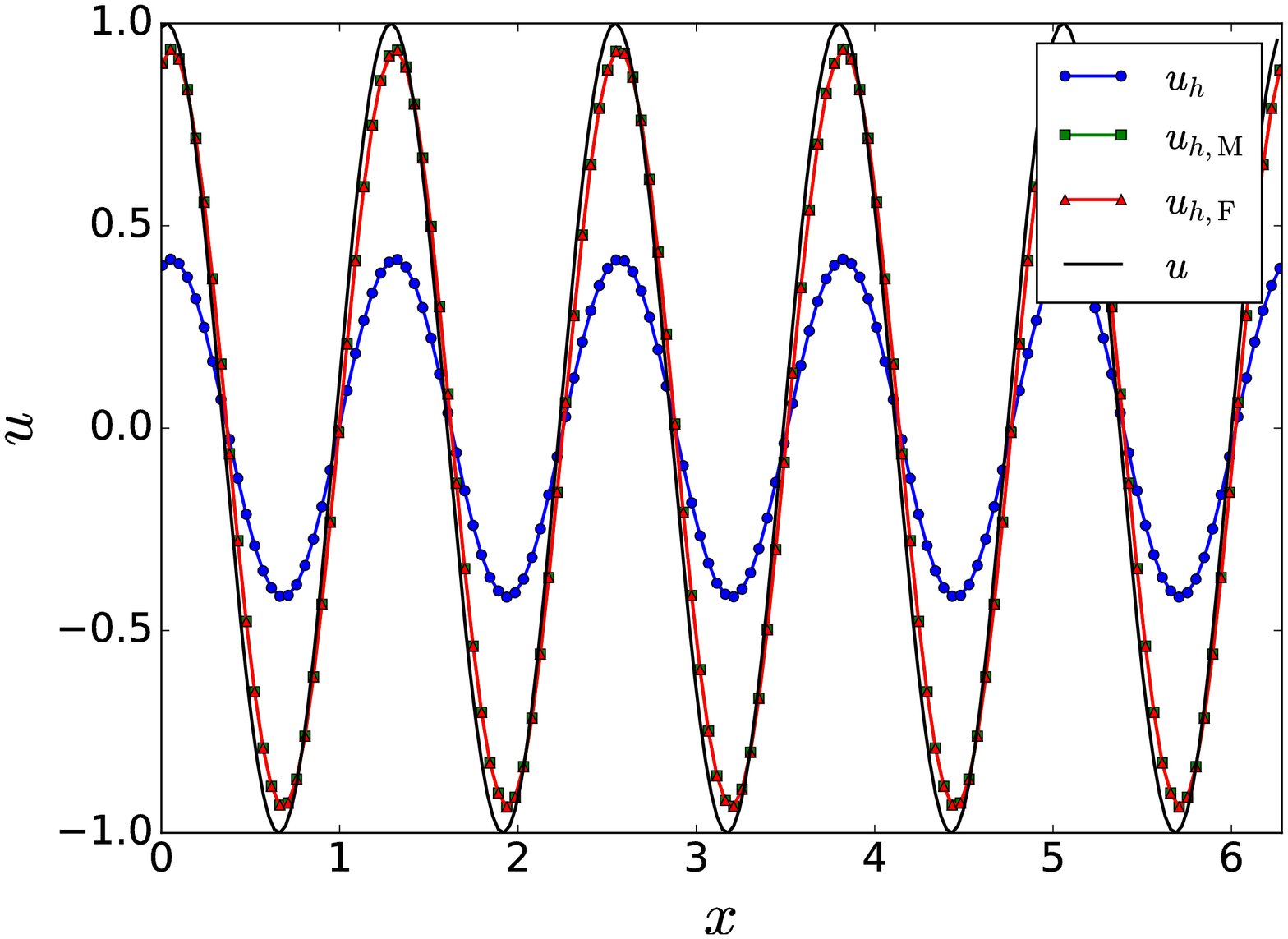}
		\caption{$u$ at $T = 1000$, $P^3$-DG, upwind flux.}
	\end{subfigure}
	\begin{subfigure}{0.45\textwidth}
		\includegraphics[width=\textwidth]{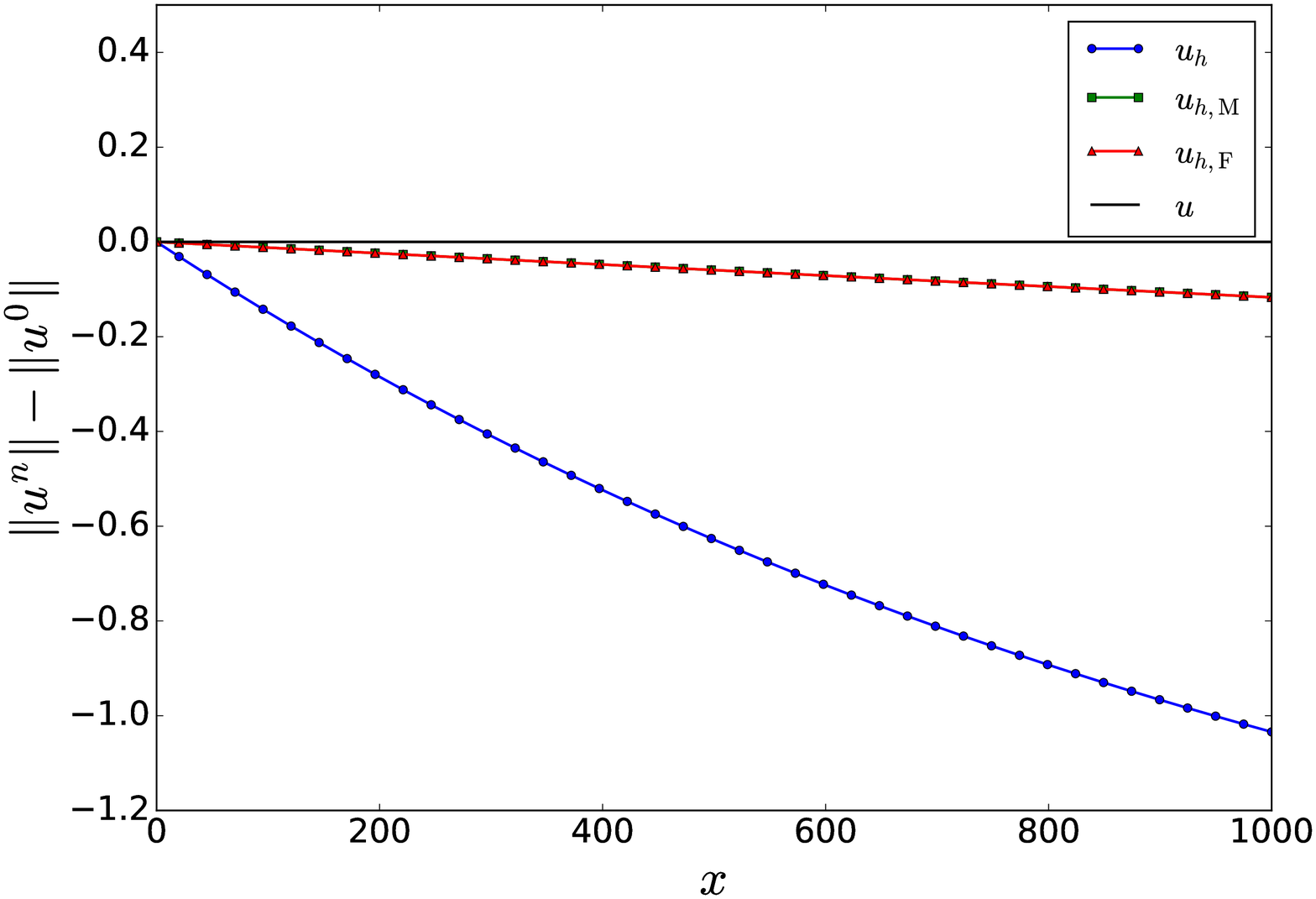}
		\caption{$\|u^n\|-\|u^0\|$.}
	\end{subfigure}
	\caption{Solution profiles and the relative $L^2$ norm. $p = k = 3$, $\alpha = -1$, $(\mu,\nu) = (0,\frac{1}{24})$, $\tau = 0.1h$, $h = \frac{2\pi}{20}$ and $u(x,0) = \sin (5x)$.}\label{fig-k3k3-energy}
\end{figure}

Besides enforcing norm decay with superviscosity, it is also possible to 
compensate the solution with anti-superviscosity to avoid the over-damping of the solution after long time simulation, which is especially notable for high frequency waves.  If the anti-superviscosity term does
not exceed the bound in Section \ref{sec-mf}, the method would still be strongly stable. 
In Figure \ref{fig-k3k3-energy}, we use $u(x,0) = \sin(5x)$ as the initial condition and use the third order time integrator with $P^3$ upwind DG method for simulation. $20$ mesh cells are used over the domain. We take $\tau = 0.1h$ and  compute to $T = 1000$. As we can see, the solution with the original third order RK method has been severely damped. While after applying anti-superviscosity, the solution maintains the amplitude reasonably well. We use a similar setting and compute with the DG scheme using the central flux. Similar difference in the solution amplitude is observed after very long time $T = 10000$. One can also see from the figure, the modified and filtering methods have a smaller decay in $L^2$ norm, and can capture the energy evolution better. However, since they are still strongly stable methods and the solution norm are decreasing with respect to time, if we take a larger final time $T$, the solution amplitude will still be damped.

\begin{figure}[h!]
	\centering
	\begin{subfigure}{0.45\textwidth}
		\includegraphics[width=\textwidth]{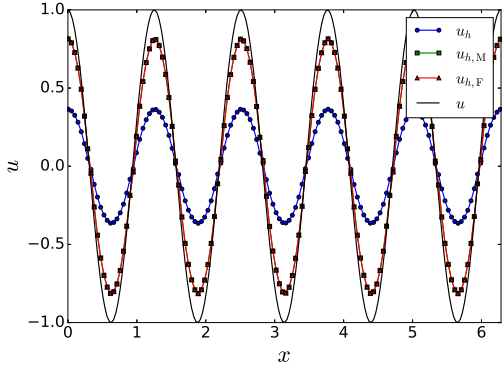}
		\caption{$u$ at $T = 10000$, $P^2$-DG, central flux.}
	\end{subfigure}
	\begin{subfigure}{0.45\textwidth}
		\includegraphics[width=\textwidth]{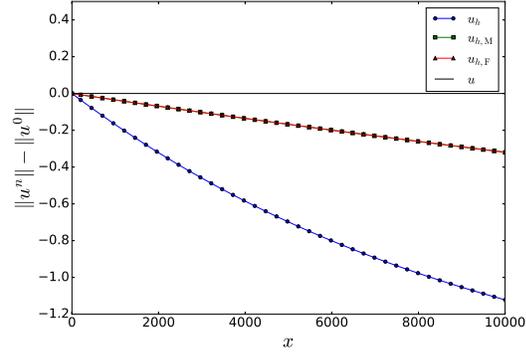}
		\caption{$\|u^n\|-\|u^0\|$.}
	\end{subfigure}
	\caption{Solution profiles and the relative $L^2$ norm. $p = 3$, $k = 2$, $\alpha = 0$, $(\mu,\nu) = (0,\frac{1}{30})$, $\tau = 0.1h$, $h = \frac{2\pi}{20}$ and $u(x,0) = \sin (5x)$.}
\end{figure}

\subsubsection{Discontinuous initial data}
\begin{figure}
	\centering
	\begin{subfigure}{0.45\textwidth}
		\includegraphics[width=\textwidth]{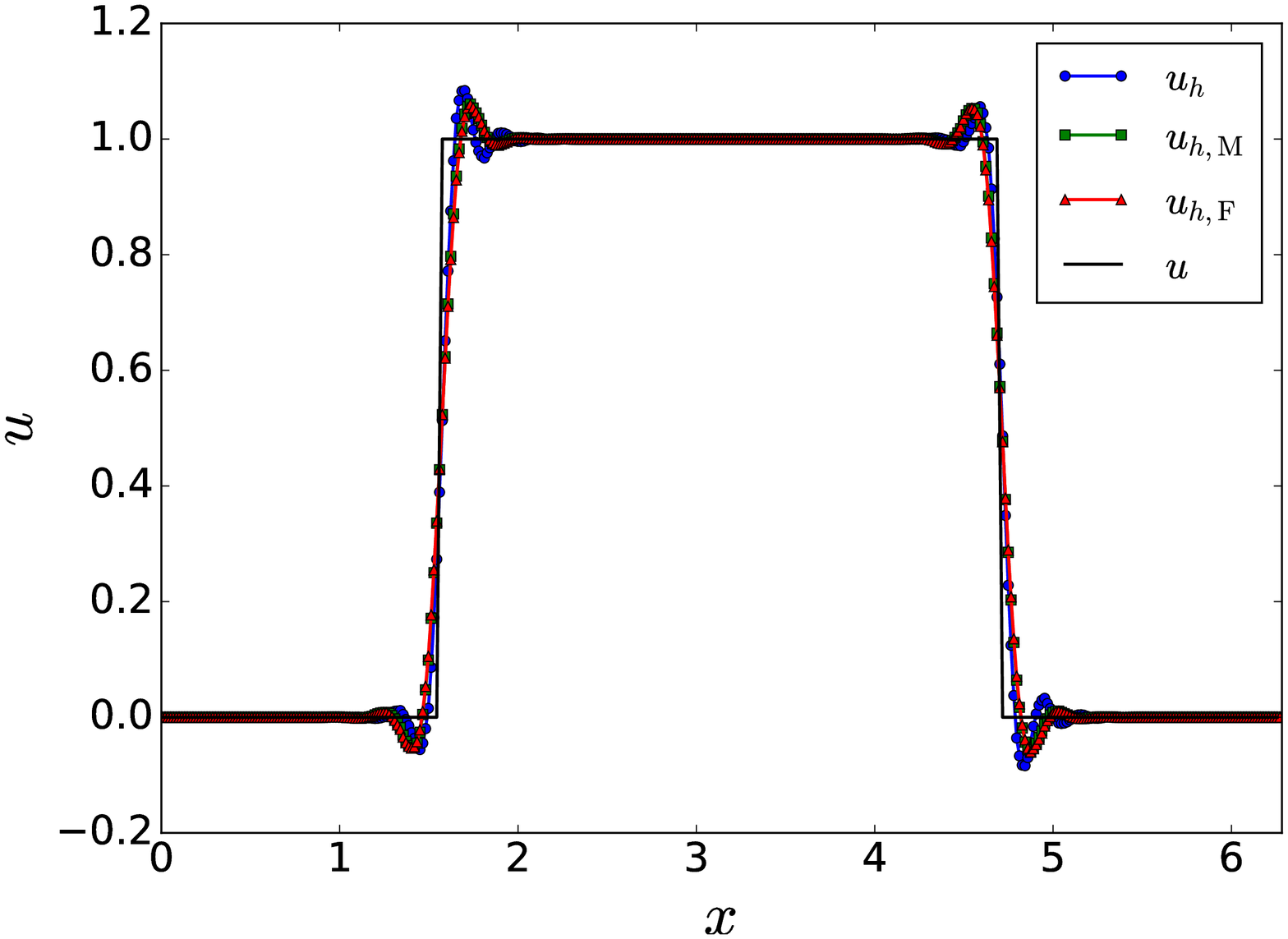}
		\caption{$u$ at $T = 2\pi$, $P^2$-DG, upwind flux.}
	\end{subfigure}
	\begin{subfigure}{0.45\textwidth}
		\includegraphics[width=\textwidth]{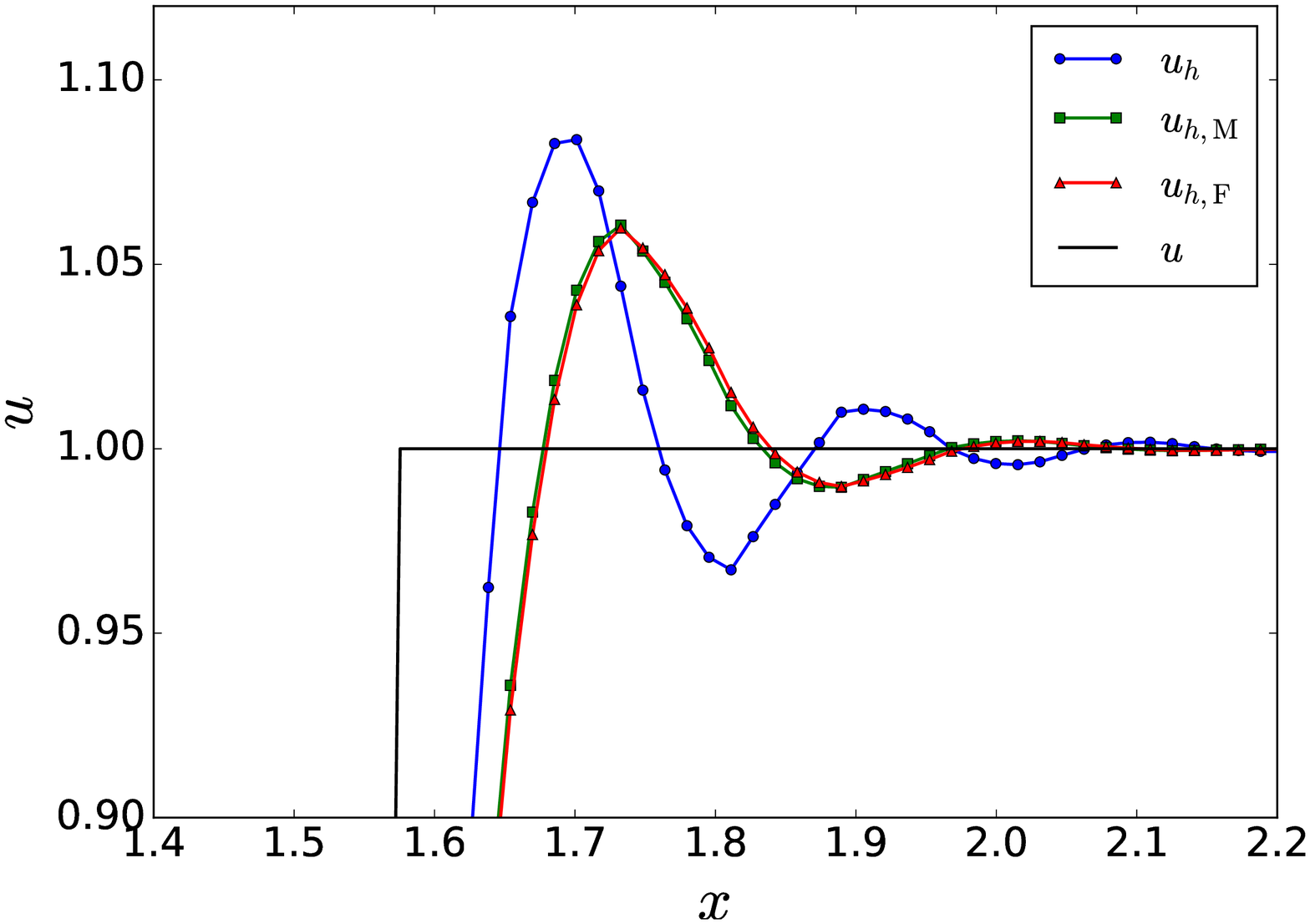}
		\caption{Zoomed in plot near the overshoot.}
	\end{subfigure}
	\caption{Solution profiles of \eqref{eq-adv} with discontinuous initial input $u(x,0) = 1_{[-\frac{\pi}{2},\frac{3\pi}{2}]}(x)$ after 1 period. $p = 3$, $k = 2$, $\alpha = -1$, $(\mu,\nu) = (0,-10)$, $\tau = 0.05h$ and $h = \frac{2\pi}{80}$. }\label{fig-disc-k2p3-80}
\end{figure}

\begin{figure}
	\centering
	\begin{subfigure}{0.45\textwidth}
		\includegraphics[width=\textwidth]{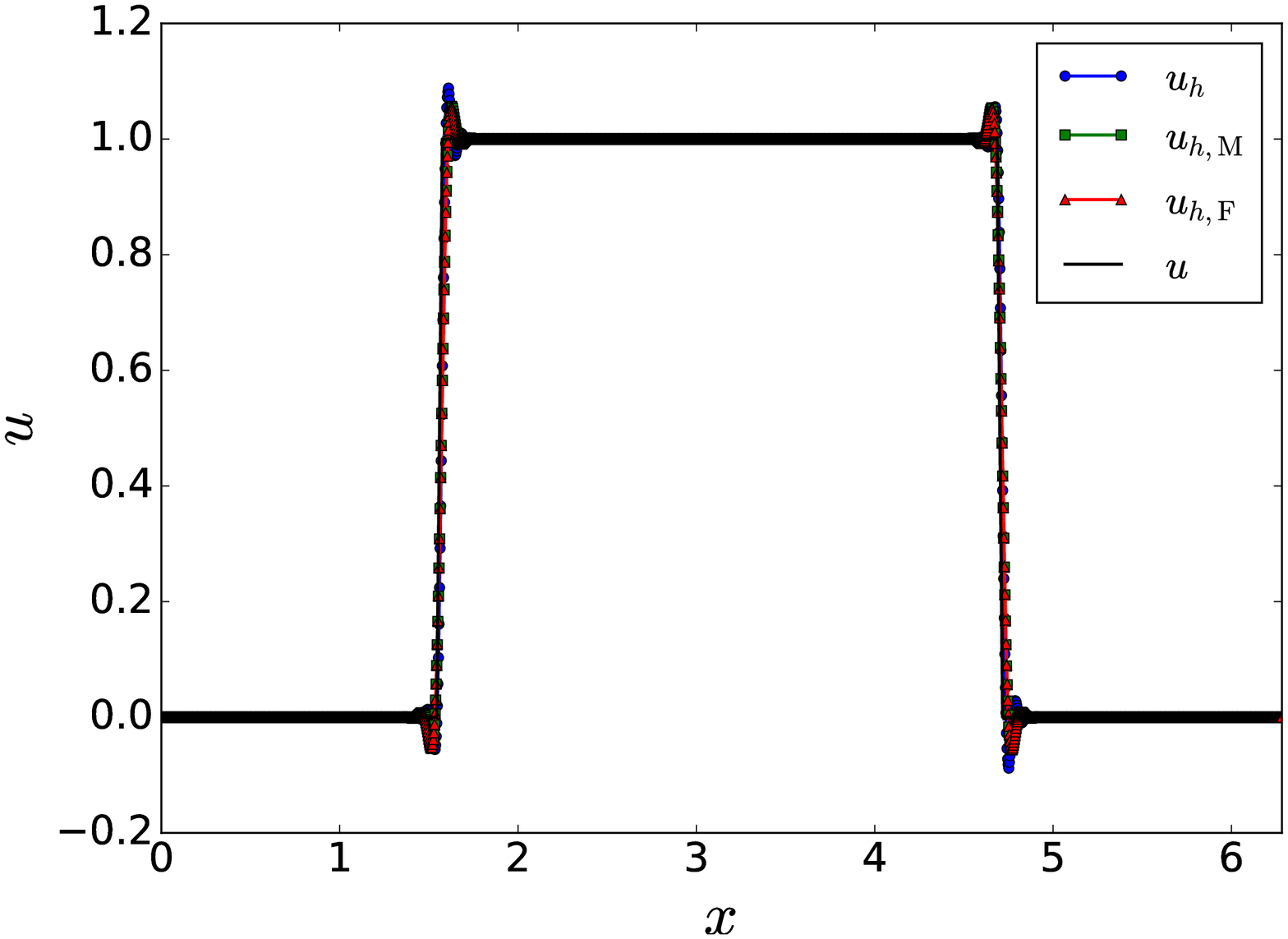}
		\caption{$u$ at $T = 2\pi$, $P^2$-DG, upwind flux.}
	\end{subfigure}
	\begin{subfigure}{0.45\textwidth}
		\includegraphics[width=\textwidth]{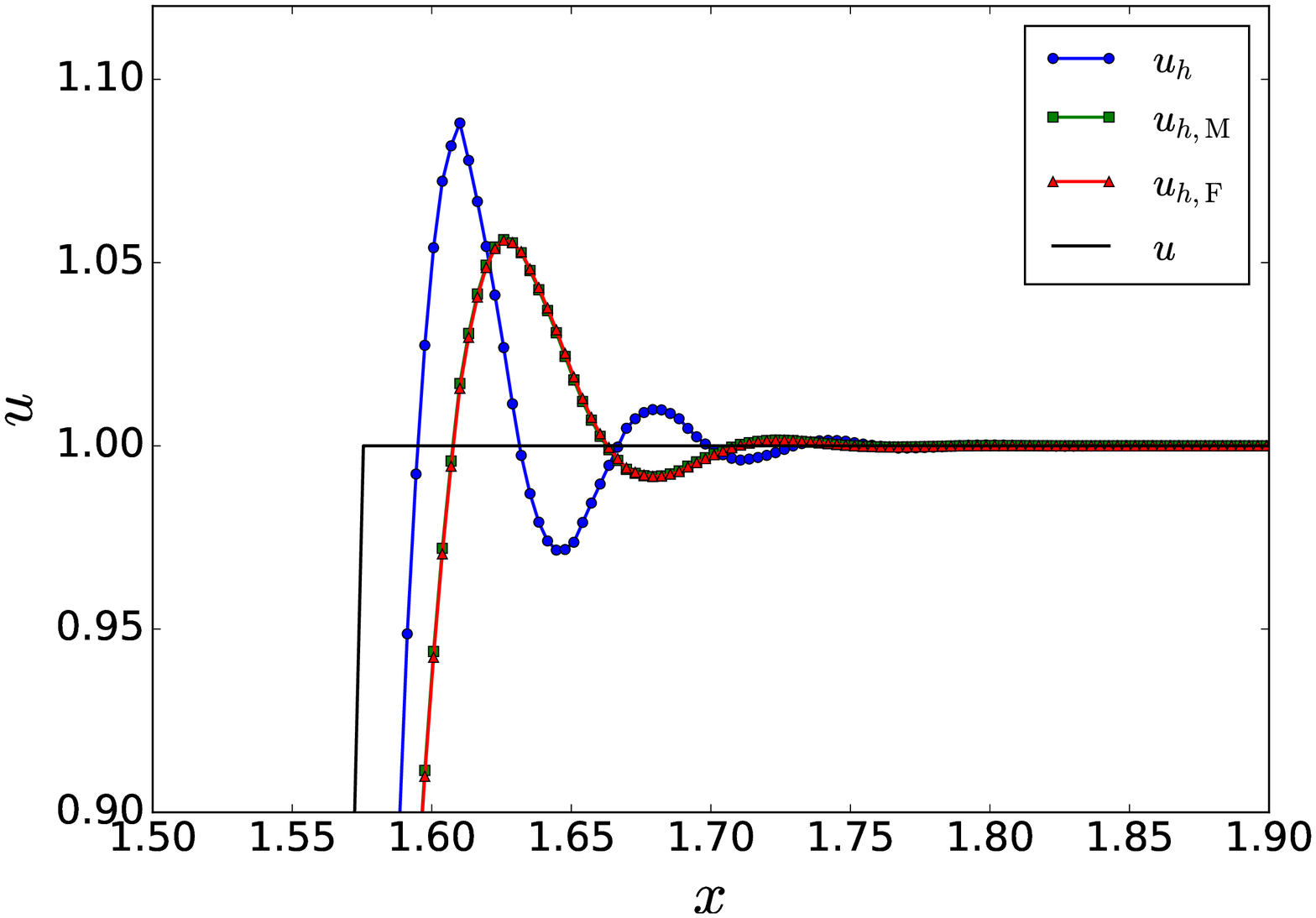}
		\caption{Zoomed in plot near the overshoot.}
	\end{subfigure}
	\caption{Solution profiles of \eqref{eq-adv} with discontinuous initial input $u(x,0) = 1_{[-\frac{\pi}{2},\frac{3\pi}{2}]}(x)$ after 1 period. $p = 3$, $k = 2$, $\alpha = -1$, $(\mu,\nu) = (0,-10)$, $\tau = 0.05h$ and $h = \frac{2\pi}{320}$. }\label{fig-disc-k2p3-320}
\end{figure}

\begin{figure}[h!]
	\centering
	\begin{subfigure}{0.32\textwidth}
		\includegraphics[width=\textwidth]{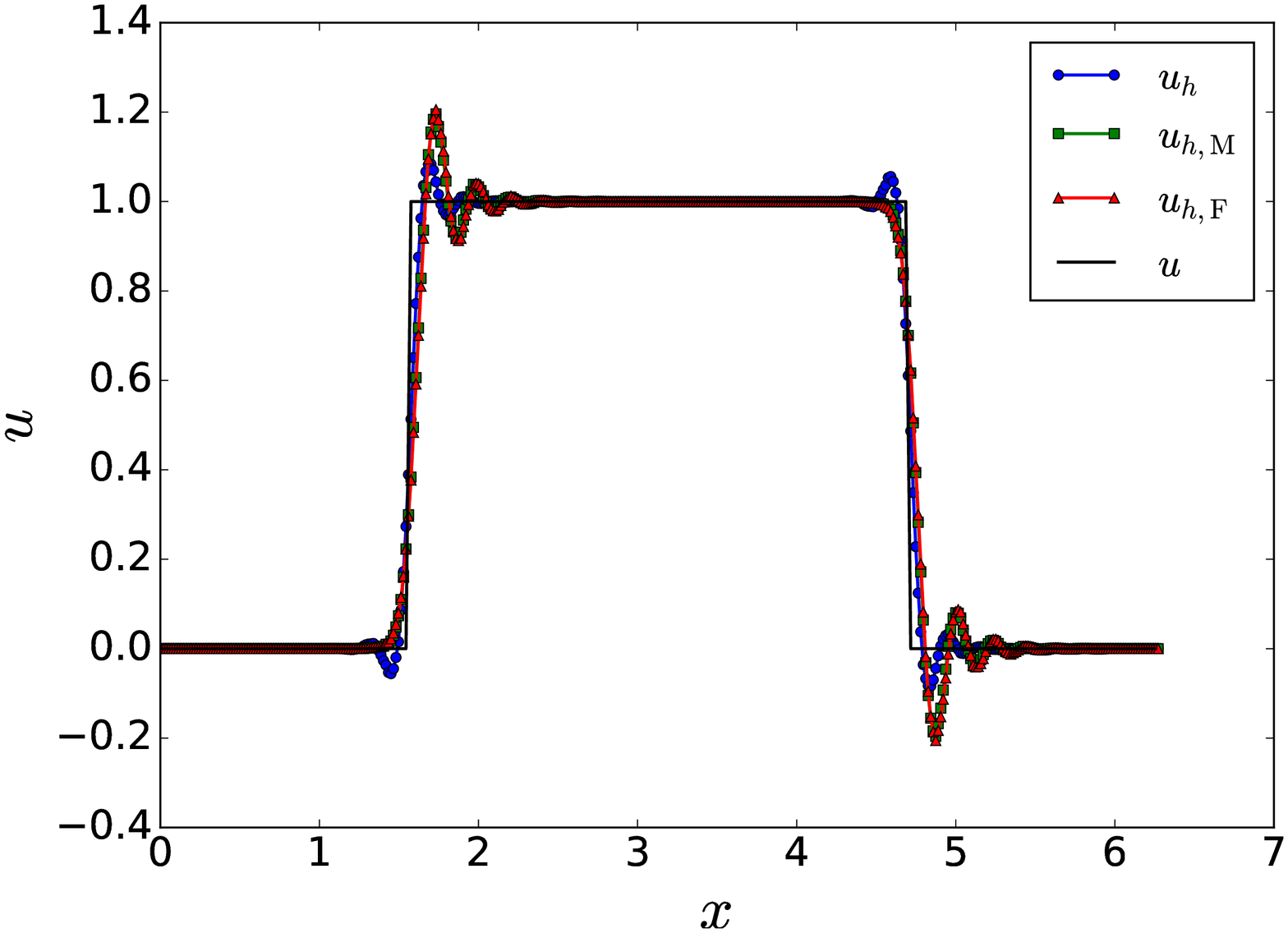}
		\caption{Solution profiles.}	
	\end{subfigure}
	\begin{subfigure}{0.32\textwidth}
		\includegraphics[width=\textwidth]{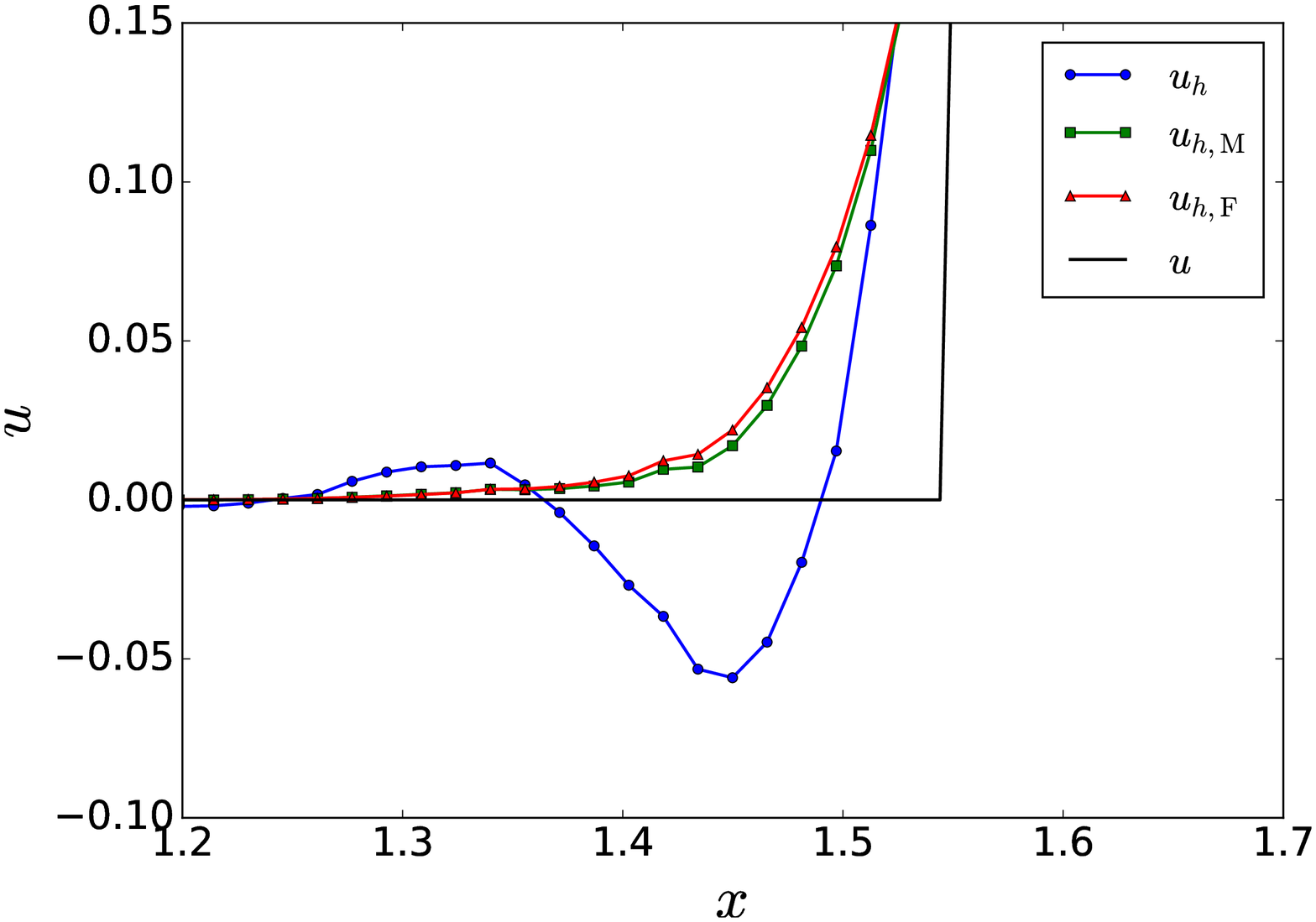}
		\caption{Undershoots.}	
	\end{subfigure}
	\begin{subfigure}{0.32\textwidth}
		\includegraphics[width=\textwidth]{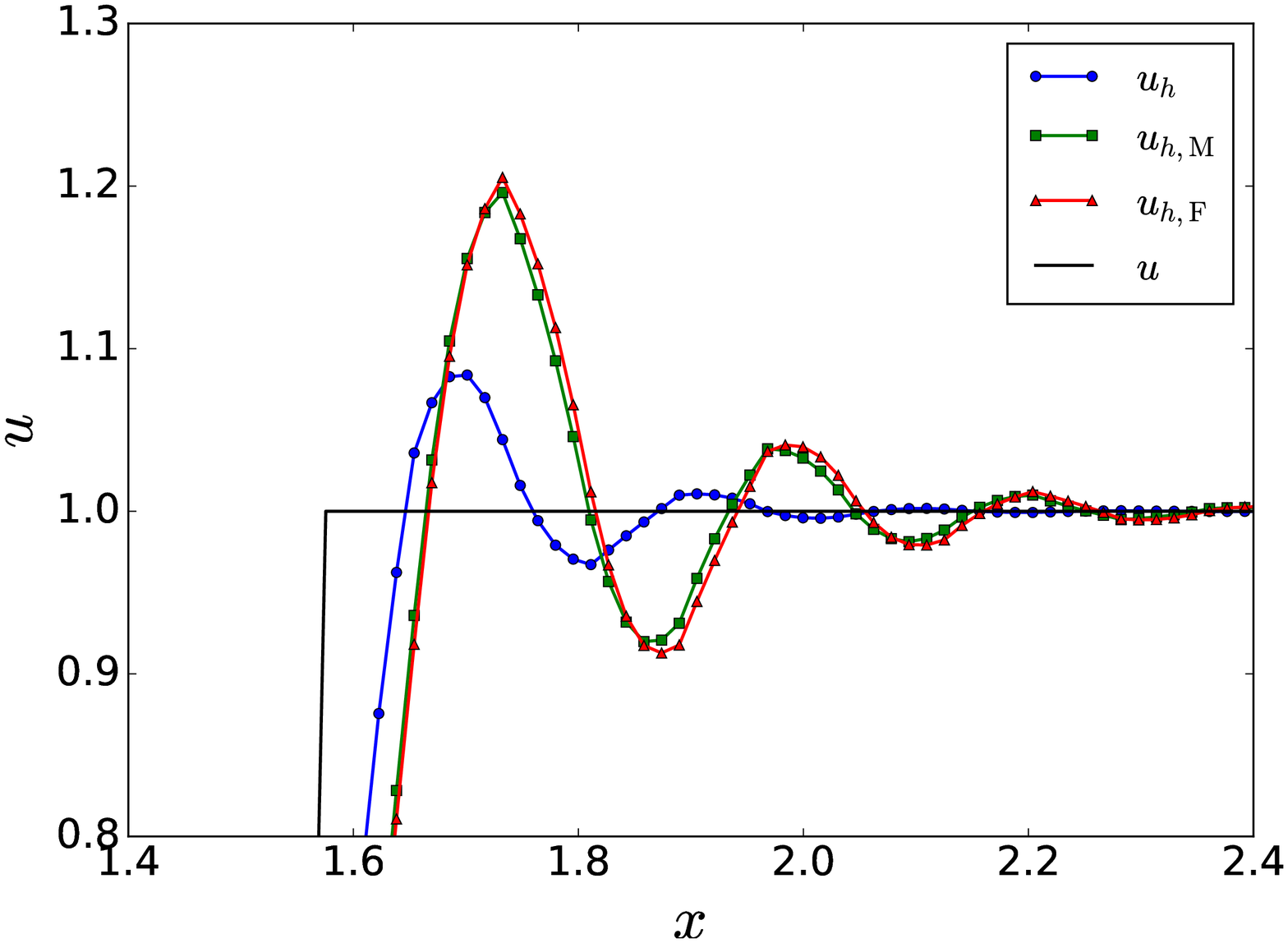}
		\caption{Overshoots.}	
	\end{subfigure}
	\caption{Solution profiles with discontinuous initial input $u(x,0) = 1_{[-\frac{\pi}{2},\frac{3\pi}{2}]}(x)$ after 1 period. $p = 3$, $k = 2$, $\alpha = -1$, $(\mu,\nu) = (1,0)$, $\tau = 0.05h$ and $h = \frac{2\pi}{80}$.}\label{fig-disper}
\end{figure}

Next, we apply the modified and filtering RK method to \eqref{eq-adv} with a discontinuous initial condition $u(x,0) = 1_{[\frac{\pi}{2},\frac{3\pi}{2}]}(x)$. When upwind DG method is used, the numerical solution may have spurious oscillations near the discontinuity. In Figure \ref{fig-disc-k2p3-80} and Figure \ref{fig-disc-k2p3-320}, we use the third order time integrator and $P^2$ upwind DG method with $80$ and $320$ mesh cells for the
simulation. Compared with those of the original methods, overshoots and undershoots of modified and filtering methods are indeed smaller, however
it seems difficult to get rid of spurious oscillations completely, since the superviscosity is only of the order $\mathcal{O}(\tau^4)$ and its effect on the final solution is quite limited. A numerical test with dispersive superviscosity is also given in Figure \ref{fig-disper}. The modified and filtered methods disperse the oscillations into one side of the discontinuity,
while the solution overshoots and undershoots do not seem to be damped. 

\begin{figure}[h!]
	\centering	
	\begin{subfigure}{0.45\textwidth}
		\includegraphics[width=\textwidth]{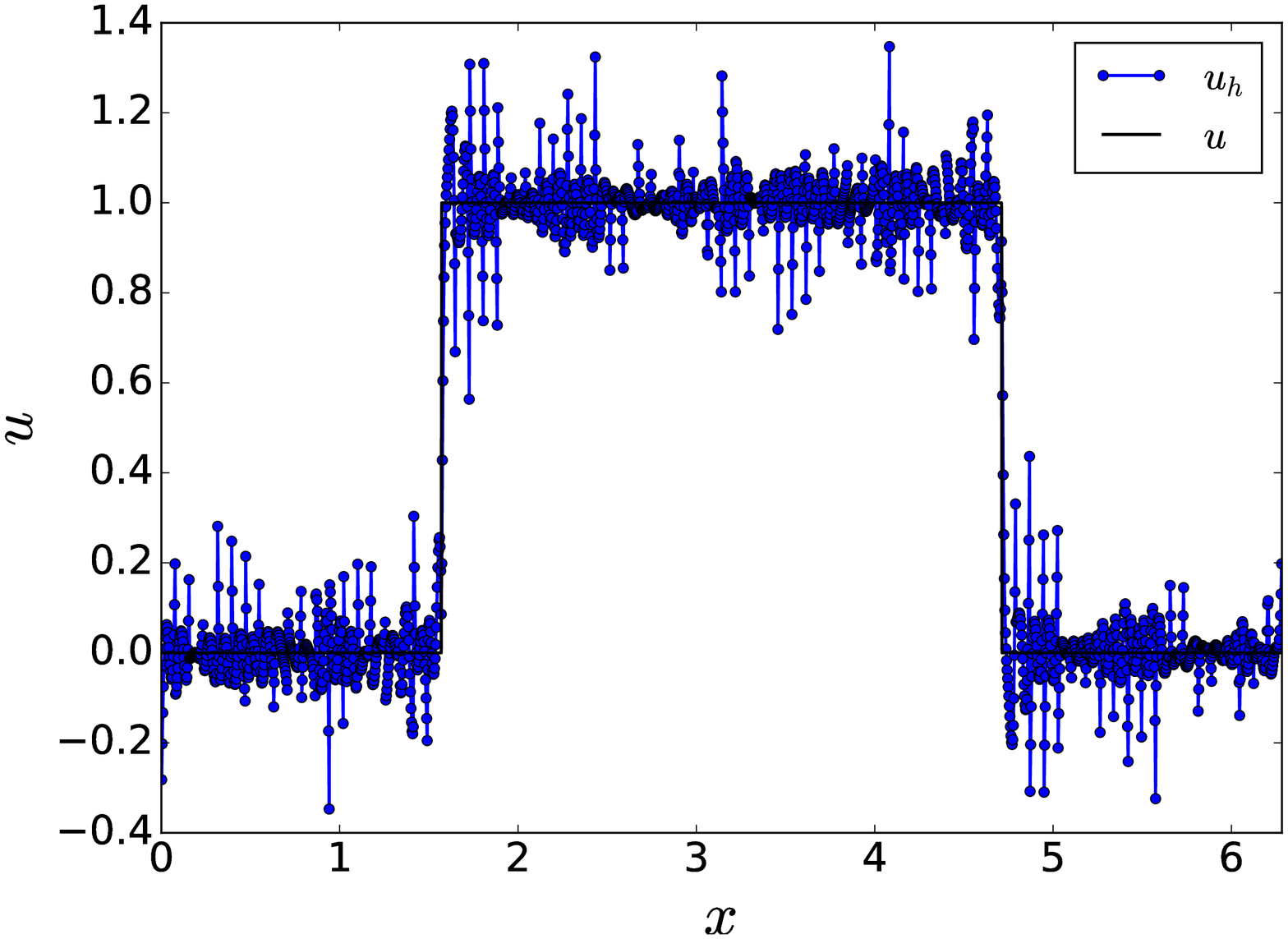}
		\caption{$u_h$, $(\mu,\nu) = (0,0)$.}
	\end{subfigure}
	\begin{subfigure}{0.45\textwidth}
	\includegraphics[width=\textwidth]{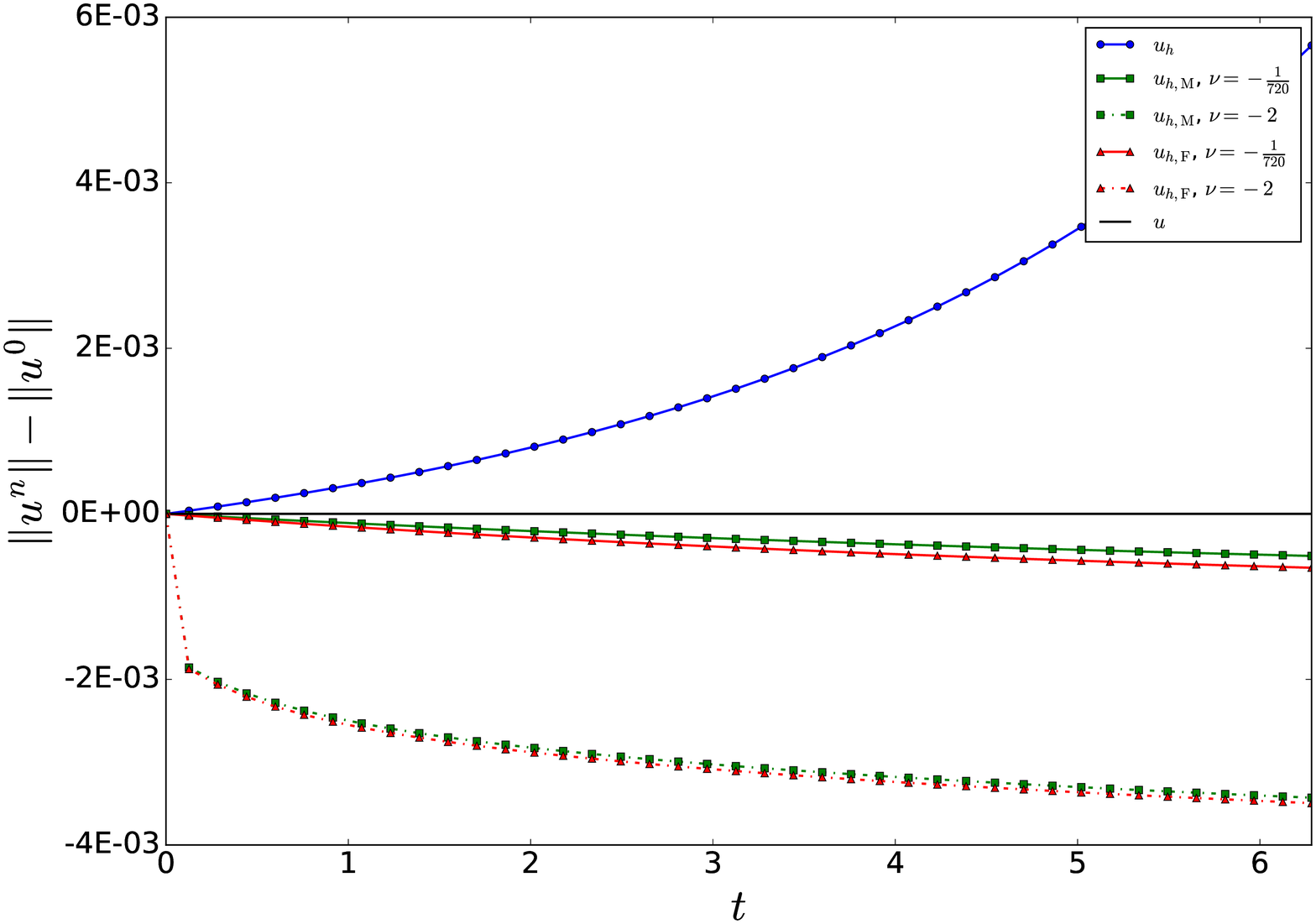}
	\caption{$\|u^n \|-\|u^0\|$.}
	\end{subfigure}
	\begin{subfigure}{0.45\textwidth}
		\includegraphics[width=\textwidth]{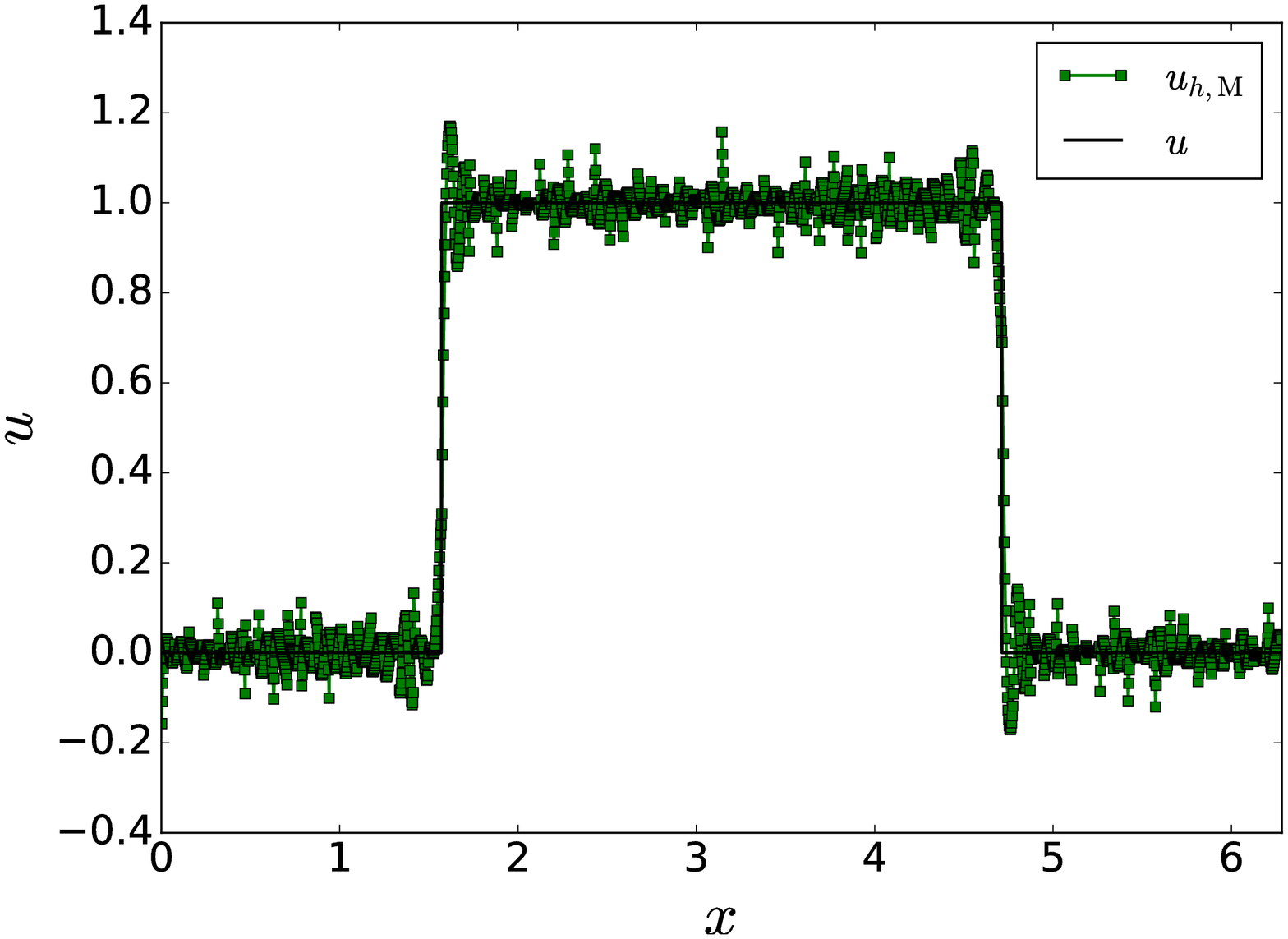}
		\caption{$u_{h,\rm{M}}$, $(\mu,\nu) = (0,-\frac{1}{720})$.}
	\end{subfigure}
	\begin{subfigure}{0.45\textwidth}
		\includegraphics[width=\textwidth]{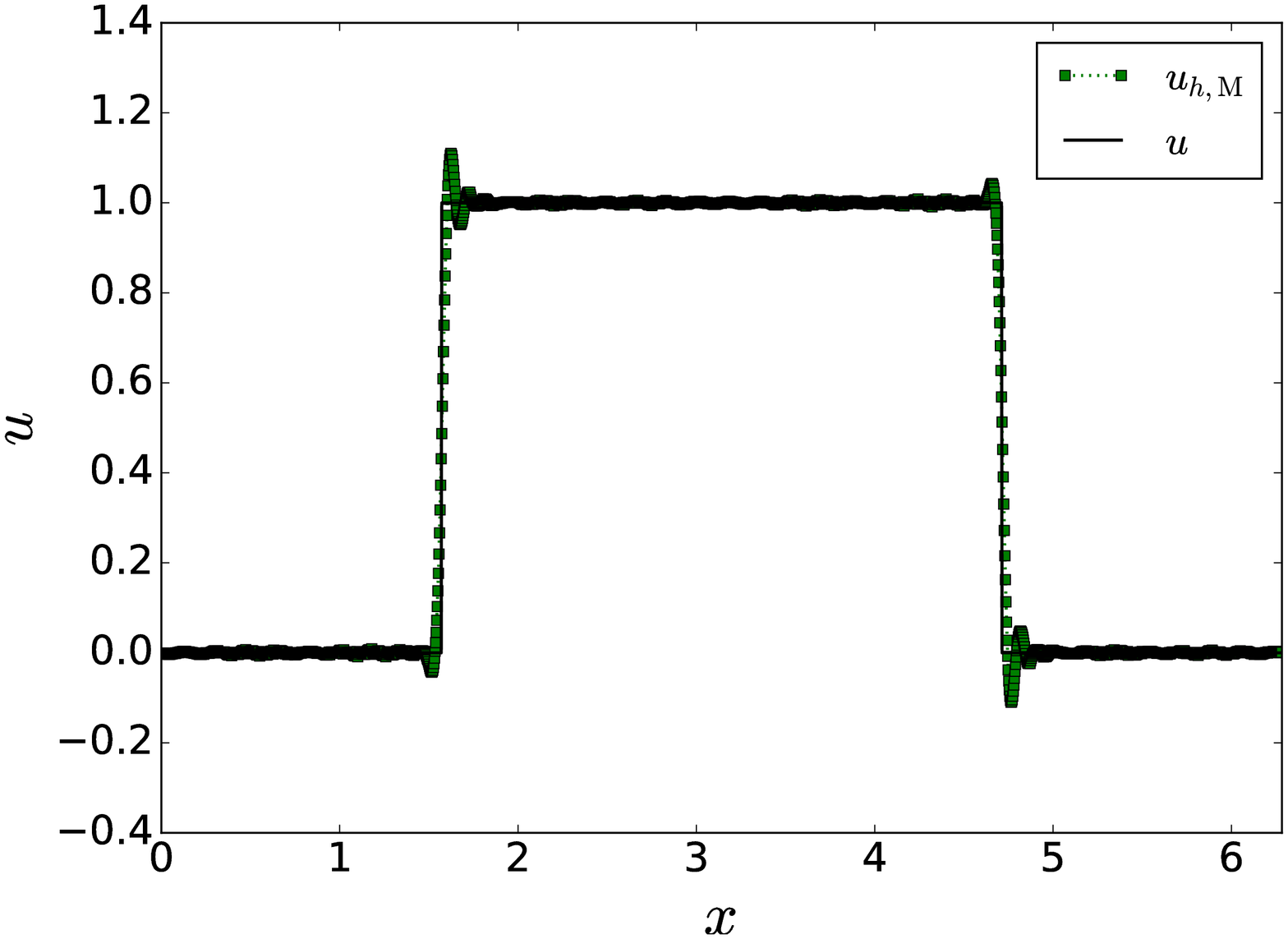}
		\caption{$u_{h,\rm{M}}$, $(\mu,\nu) = (0,-2)$.}
	\end{subfigure}
	\begin{subfigure}{0.45\textwidth}
		\includegraphics[width=\textwidth]{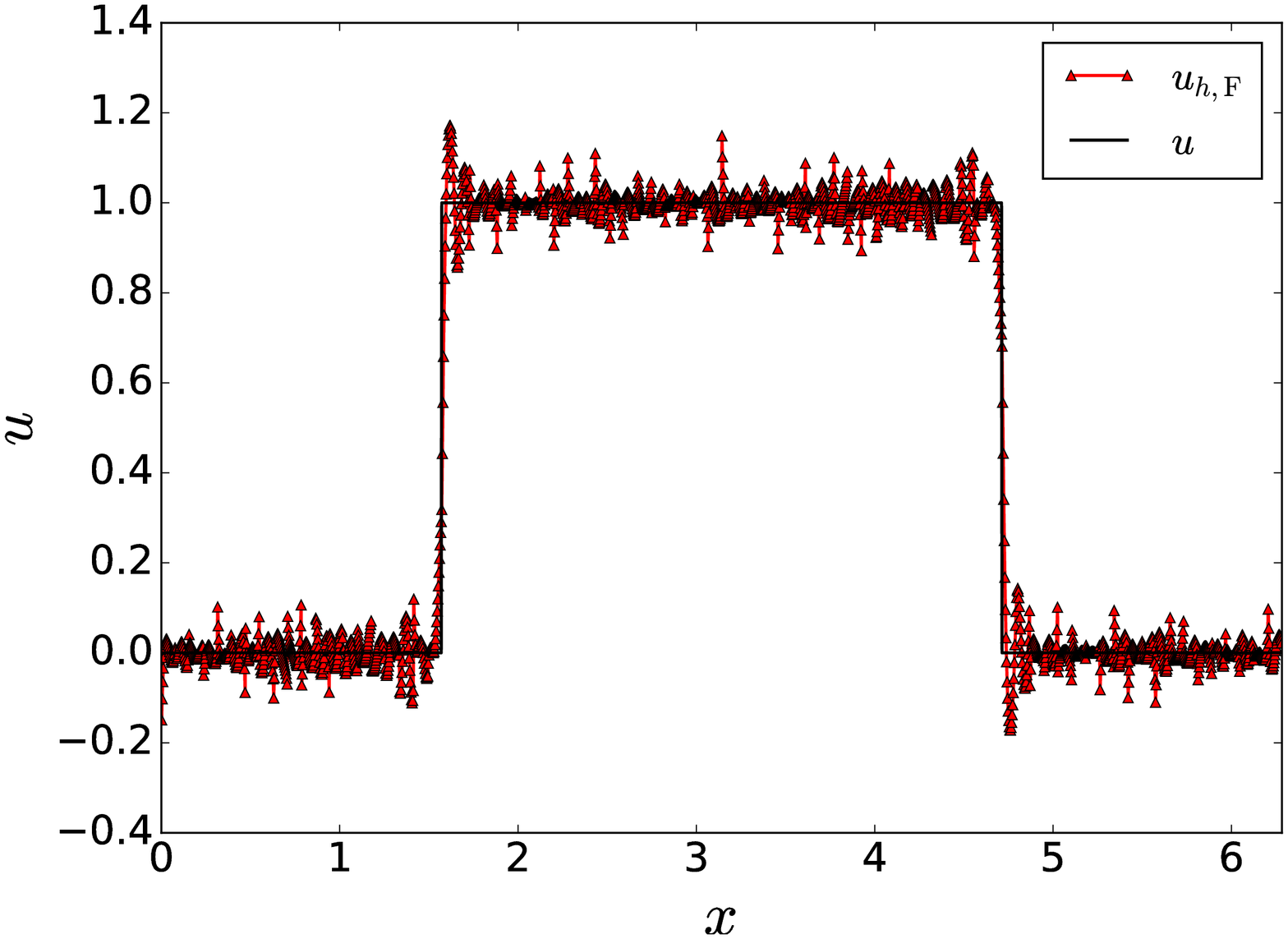}
		\caption{$u_{h,\rm{F}}$, $(\mu,\nu) = (0,-\frac{1}{720})$.}
	\end{subfigure}
	\begin{subfigure}{0.45\textwidth}
		\includegraphics[width=\textwidth]{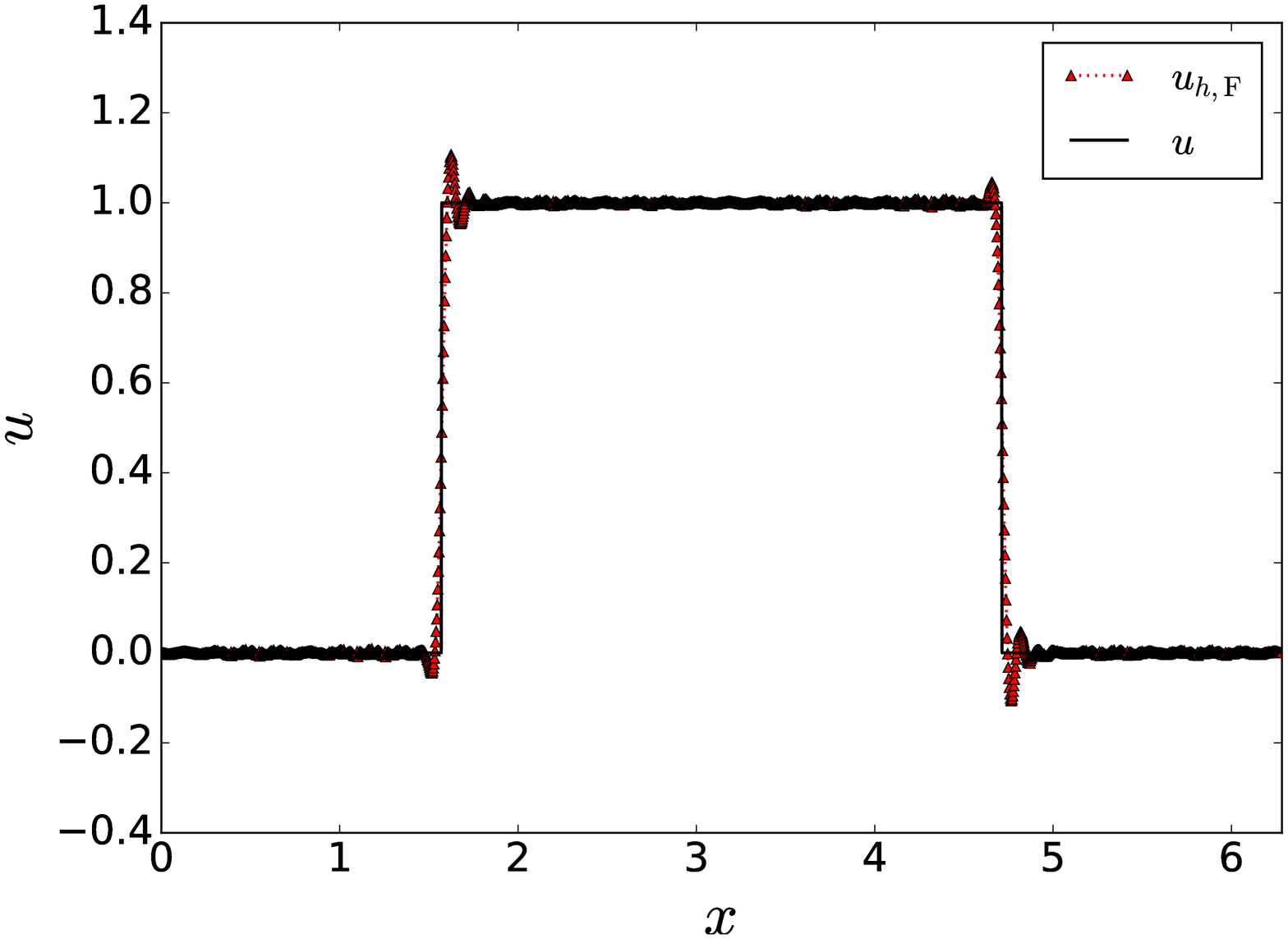}
		\caption{$u_{h,\rm{F}}$, $(\mu,\nu) = (0,-2)$.}
	\end{subfigure}
	\caption{Solution profiles with discontinuous initial input $u(x,0) = 1_{[-\frac{\pi}{2},\frac{3\pi}{2}]}(x)$ after 1 period. $ p = 3$, $k = 2$, $\alpha = 0$, $(\mu,\nu) = (0,-10)$, $\tau = 0.05h$, $h = \frac{2\pi}{80}$. }\label{fig-disc-k4p5}
\end{figure}

We then switch to the central flux to test the problem. See Figure \ref{fig-disc-k4p5}, where $P^4$ polynomials with the fifth order time integrator is used. The total number of mesh cell is $N = 80$ and the time step is $\tau = 0.05h$. With the usual time integrator the solution is highly oscillatory and the solution norm increases with respect to time. When we add the superviscosity with $\nu = -\frac{1}{720}$, the solution norm decays and the solution is slightly less oscillatory compared with the original ones. When a stronger superviscosity with $\nu = -2$ is applied to the problem, the solution seems to be smoothed out.


\subsection{DG method for Burgers equation}\label{examp-burgers}
Finally, we test the adaptive filter in Section \ref{sec-nonlinear} with $D \approx \tau^{k^*} \partial_x^{k^*}$ to the Burgers equation 
\begin{equation}
	u_t + f(u)_x = 0, \qquad f(u) = \hf u^2.
\end{equation}
The DG spatial discretization is used for spatial discretization. In other words, we look for $u_h \in V_h$, such that for all $v_h \in V_h$, 
\begin{equation}
	\int_{I_j} (u_h)_t v dx -\int_{I_j} f(u_h)v_x dx + \widehat{f}_{j+\hf} v_{j+\hf}^- - \widehat{f}_{j-\hf} v_{j-\hf}^- = 0, \qquad \forall v \in V_h.
\end{equation} 
Here we take the entropy conservative flux $\widehat{f} = \frac{1}{6}\left((u_h^+)^2 + u_h^+u_h^- + (u_h^-)^2\right)$. We will use $P^2$-DG with the second order RK method, and $P^4$-DG with the six-stage fifth order Fehlberg method  \cite{fehlberg1970klassische} for time discretization.

First we test the accuracy for the Burgers equation with adaptive superviscosity. The final time is set as $ T = 0.3$ and the time step is $\tau = 0.05h$. The asymptotic convergence rate is two for the second order method and five for the fifth order method. 

 \begin{table}[h!]
	\centering
	\begin{tabular}{c|c| c | c | c | c | c |c | c }
		\hline
		\hline
		$p$&$k$&	$N$  &$L^1$ error& order &$L^2$ error &order&$L^\infty$ error &order\\
		\hline
		\hline
		\multirow{5}{*}{$2$}&
		\multirow{5}{*}{$2$}
	    &40&6.2120E-05&   -&4.0719E-05&   -&1.0567E-04& - \\
		&&80&7.7688E-06& 3.00&5.2015E-06& 2.97&1.6568E-05& 2.67 \\
		&&160&1.2366E-06& 2.65&7.3246E-07& 2.83&2.1921E-06& 2.92 \\
		&&320&2.6101E-07& 2.24&1.3552E-07& 2.43&2.8687E-07& 2.93 \\
		&&640&7.2524E-08& 1.85&3.6789E-08& 1.88&4.0667E-08& 2.82 \\
		&&1280&1.7681E-08& 2.04&8.9232E-09& 2.04&7.8042E-09& 2.38 \\
		&&2560&4.3998E-09& 2.01&2.2188E-09& 2.01&1.9419E-09& 2.01 \\
		\hline
		\hline
		\multirow{5}{*}{$5$}&
		\multirow{5}{*}{$4$}
		&20&1.3810E-06&   -&1.1148E-06&   -&2.4683E-06&   - \\
		&&40&3.3035E-08& 5.39&2.8230E-08& 5.30&7.0330E-08& 5.13 \\
		&&80&1.0107E-09& 5.03&9.5056E-10& 4.89&4.1795E-09& 4.07 \\
		&&160&3.1150E-11& 5.02&2.9522E-11& 5.01&1.3023E-10& 5.00 \\
		\hline
		\hline
	\end{tabular}
\caption{Accuracy test in Section \ref{examp-burgers} with Burgers equation at $T = 0.3$.}
\end{table}

\begin{figure}[h!]
	\centering
	\begin{subfigure}{0.45\textwidth}
		\includegraphics[width=\textwidth]{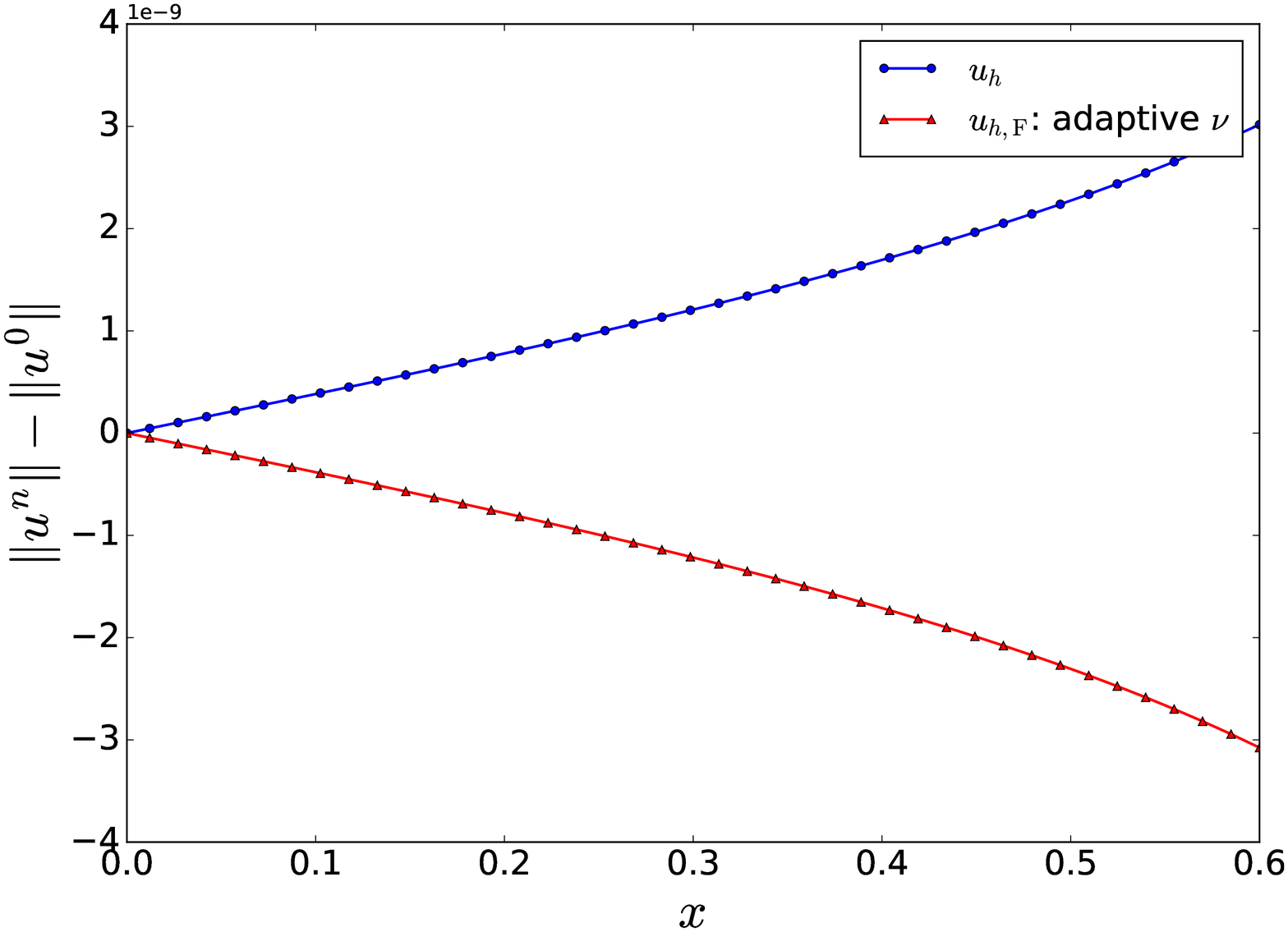}
		\caption{$p = k = 2$, $\|u^n \|-\|u^0\|$.}
	\end{subfigure}
	\begin{subfigure}{0.45\textwidth}
		\includegraphics[width=\textwidth]{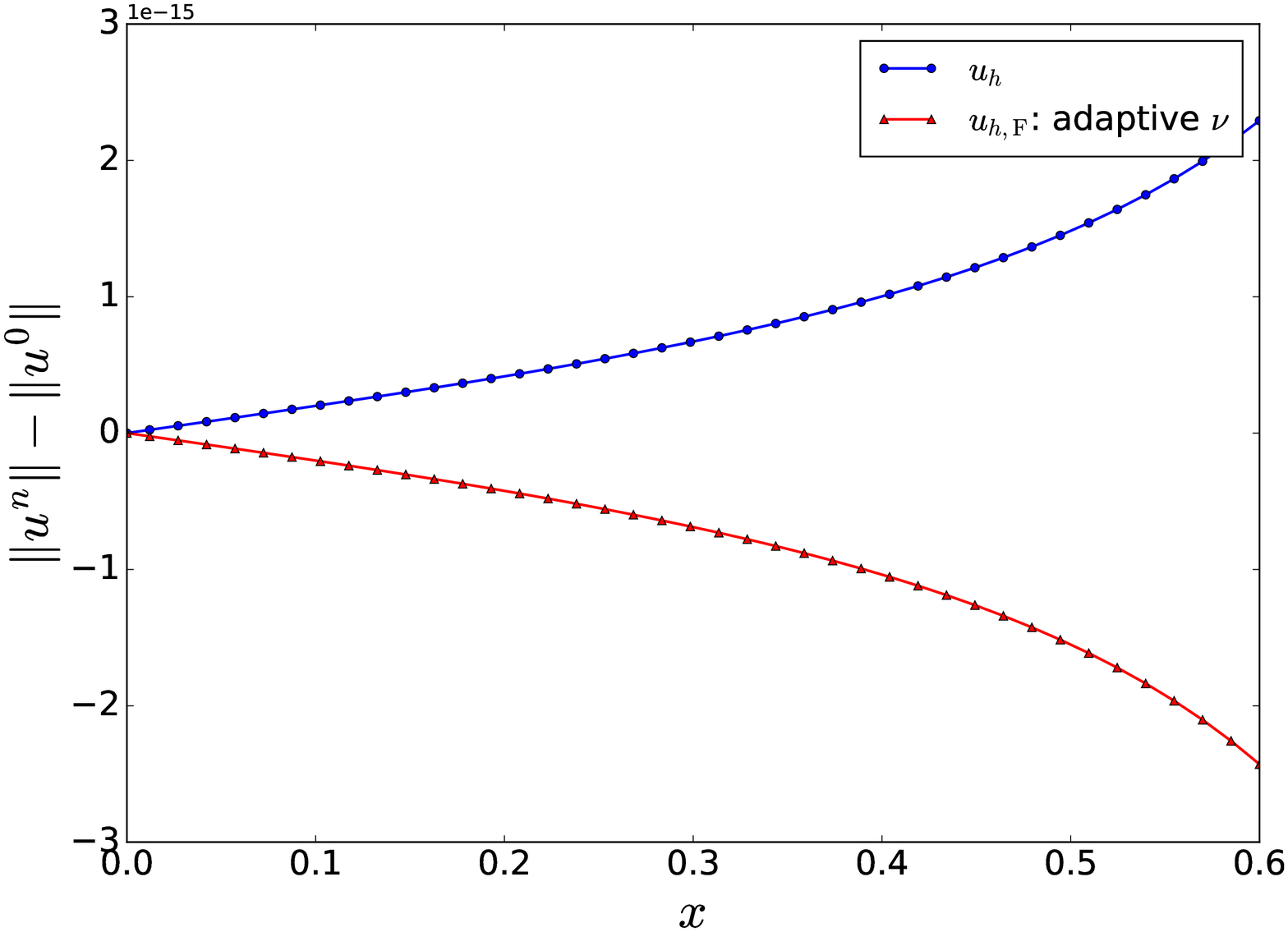}
		\caption{$p = 5$, $k = 4$, $\|u^n\|-\|u^0\|$.}
	\end{subfigure}
	\caption{Relative $L^2$ norm of solutions to Burgers equation with entropy-conserving flux. $N = 80$, $\tau = 0.05 h$.}\label{fig-adaptenergy}
\end{figure}

We then fix $N = 80$ and monitor the relative $L^2$ norm of the solution. From Figure \ref{fig-adaptenergy}, we see that for the original method, the solution norm increases with respect to time. While the solution norm decreases with the adaptive filtering method. 

In the end, we present solutions at $T = 1.5$, when the shock has been developed at $x = \pi$. The solution of the original method is highly oscillatory. Although the adaptive superviscosity damps the $L^2$ norm,
the superviscosity may not be strong enough to suppress oscillations. When a strong superviscosity term is applied through the filter, the solution seems to be smoothed out away from the discontinuity. 

\begin{figure}[h!]
	\centering
	\begin{subfigure}{0.32\textwidth}
		\includegraphics[width=\textwidth]{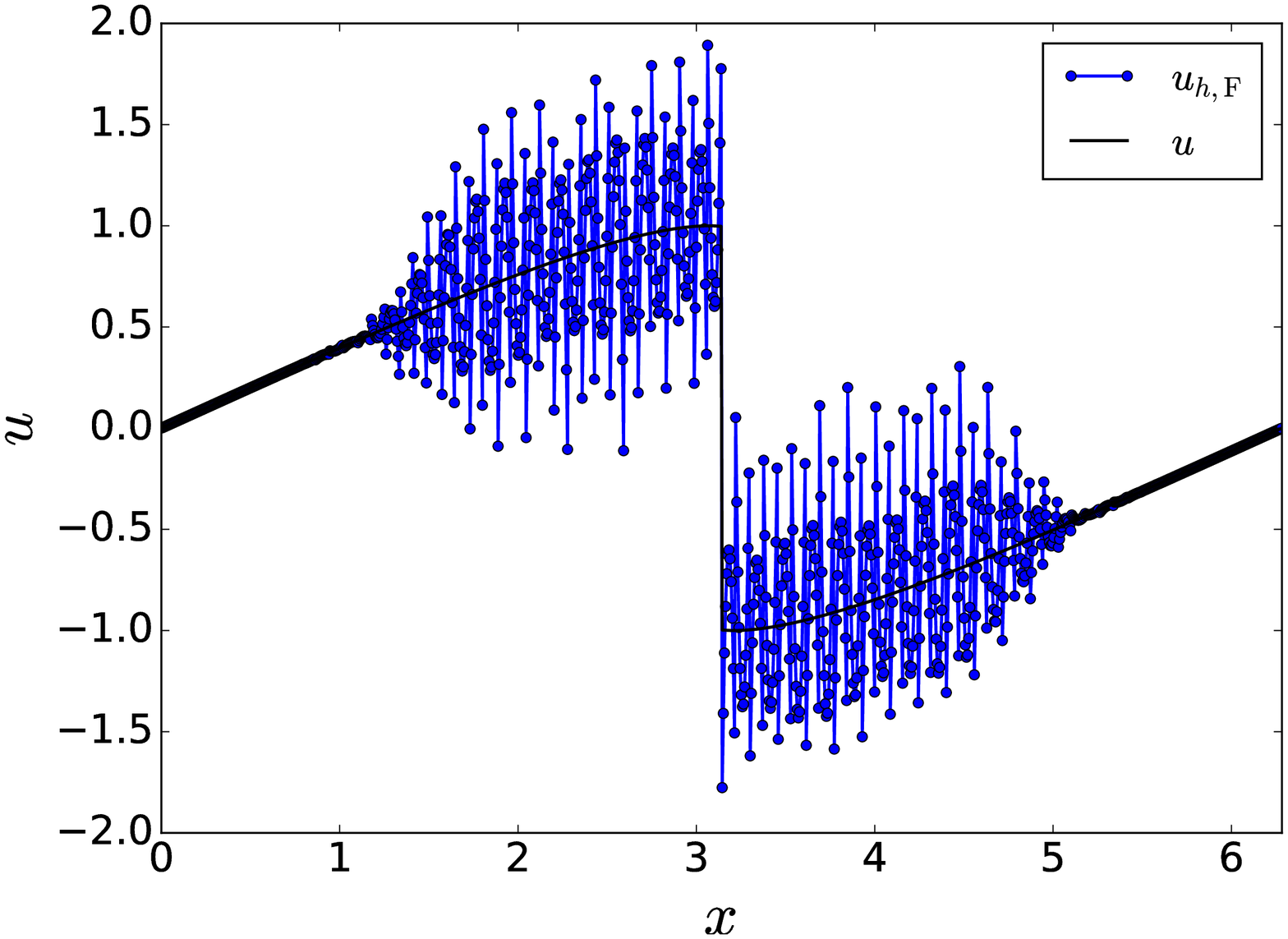}
		\caption{$\nu = 0$.}	
	\end{subfigure}
	\begin{subfigure}{0.32\textwidth}
		\includegraphics[width=\textwidth]{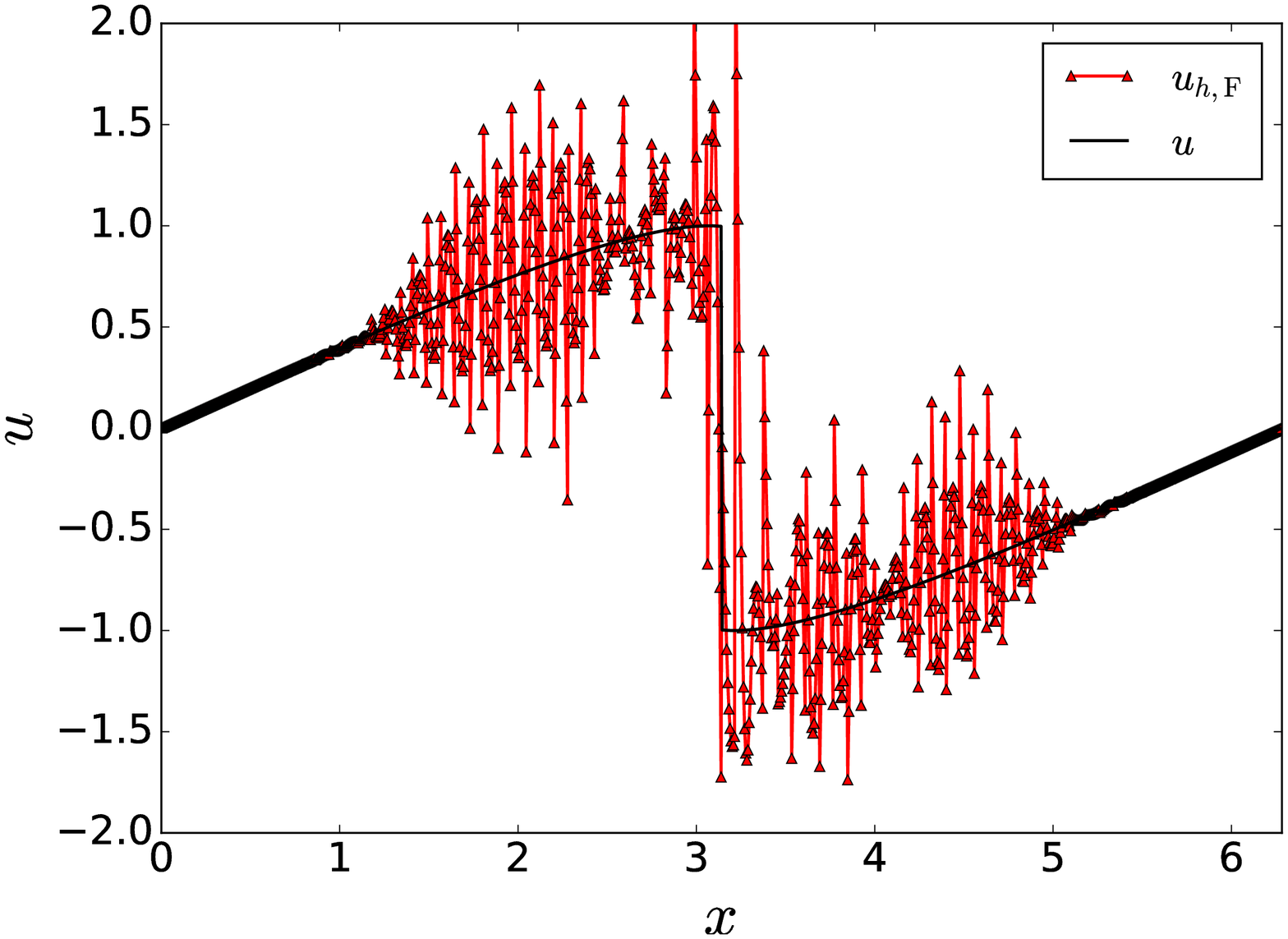}
		\caption{Adaptive $\nu$.}	
	\end{subfigure}
	\begin{subfigure}{0.32\textwidth}
		\includegraphics[width=\textwidth]{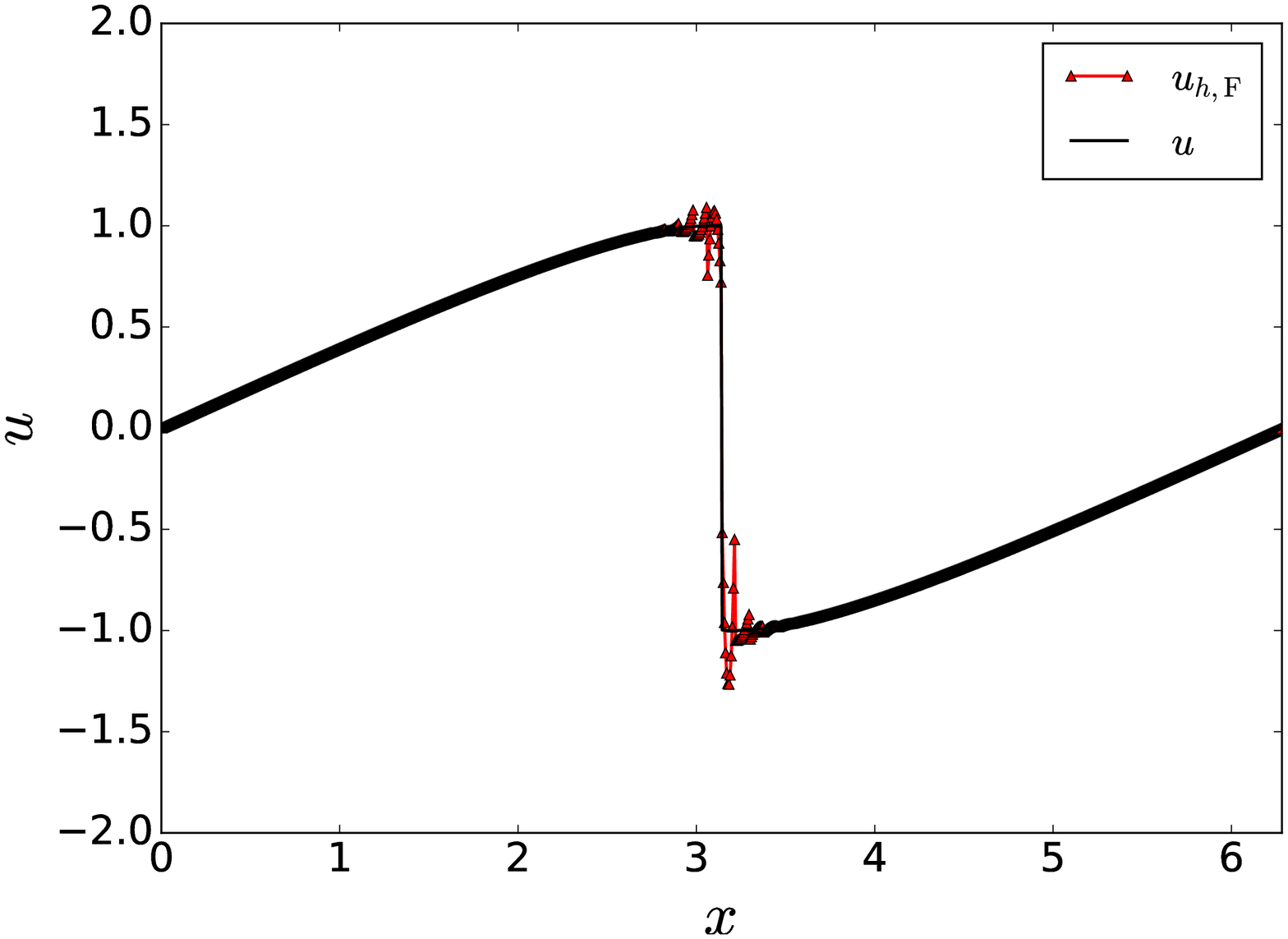}
		\caption{$\nu = -10$.}	
	\end{subfigure}
	\caption{Solution profiles at $T = 1.5$, $p = 2$, $ k = 2$, $N = 80$.}
\end{figure}
\begin{figure}[h!]
	\centering
	\begin{subfigure}{0.32\textwidth}
		\includegraphics[width=\textwidth]{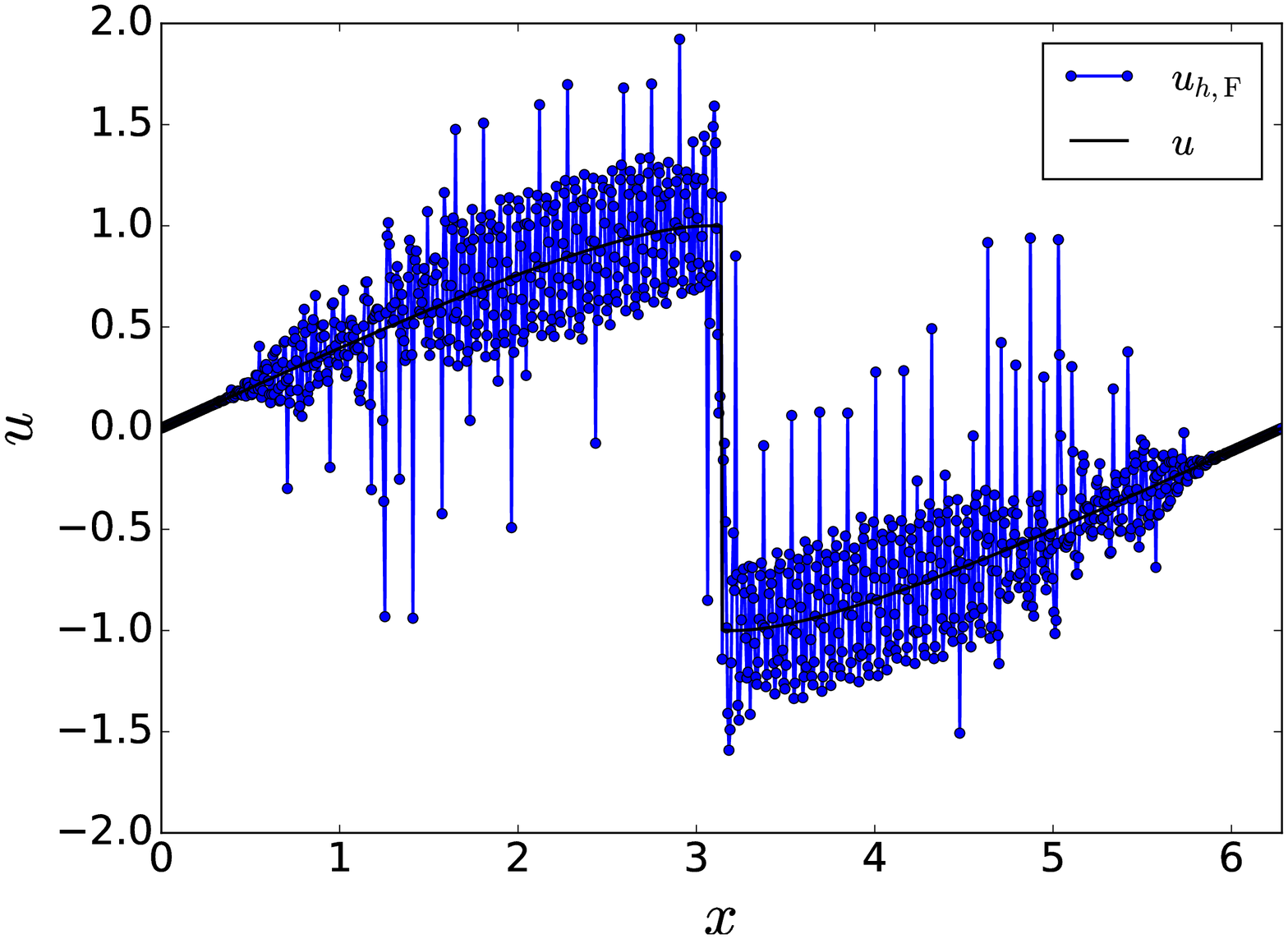}
		\caption{$\nu = 0$.}	
	\end{subfigure}
	\begin{subfigure}{0.32\textwidth}
		\includegraphics[width=\textwidth]{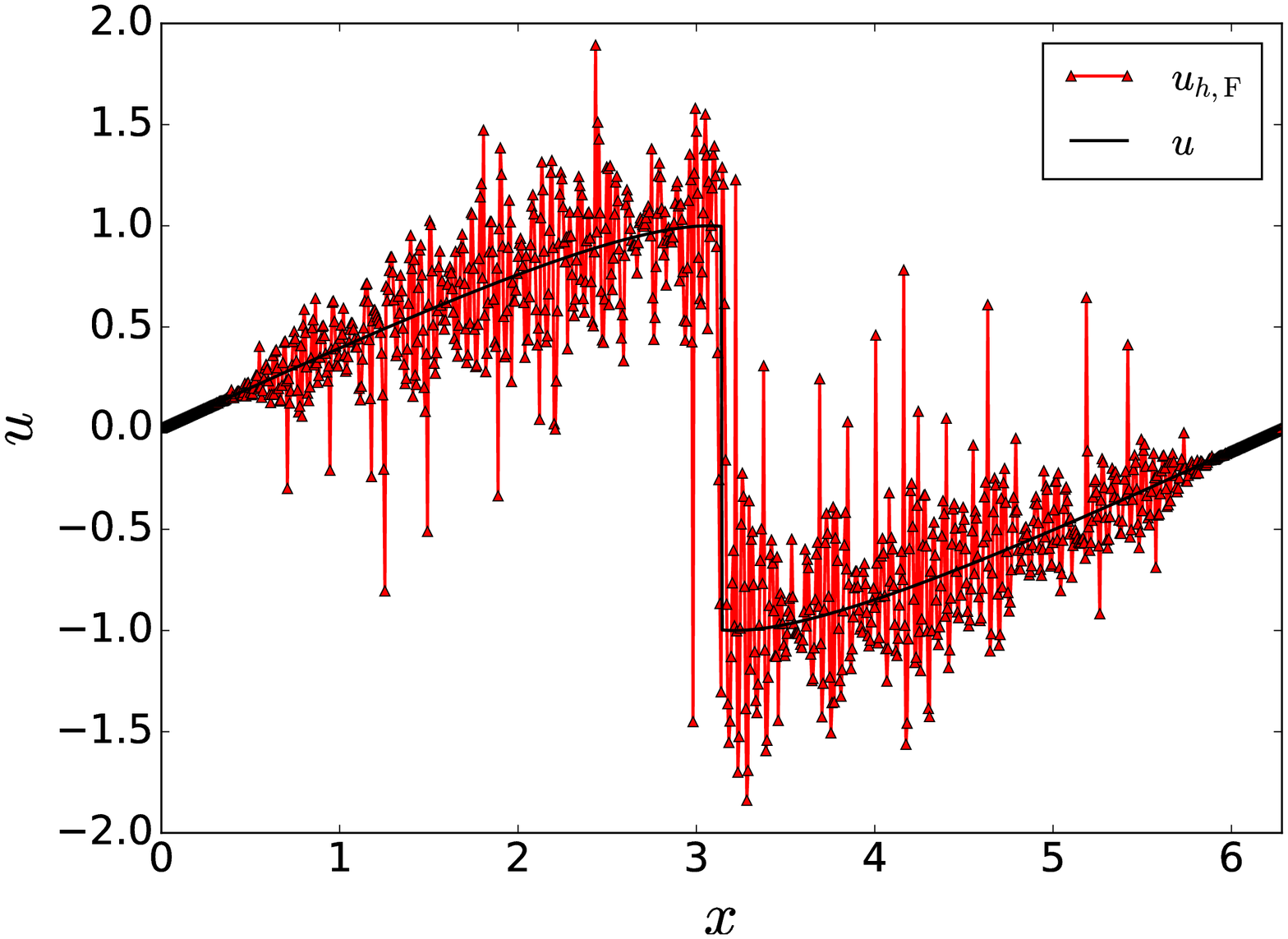}
		\caption{Adaptive $\nu$.}	
	\end{subfigure}
	\begin{subfigure}{0.32\textwidth}
		\includegraphics[width=\textwidth]{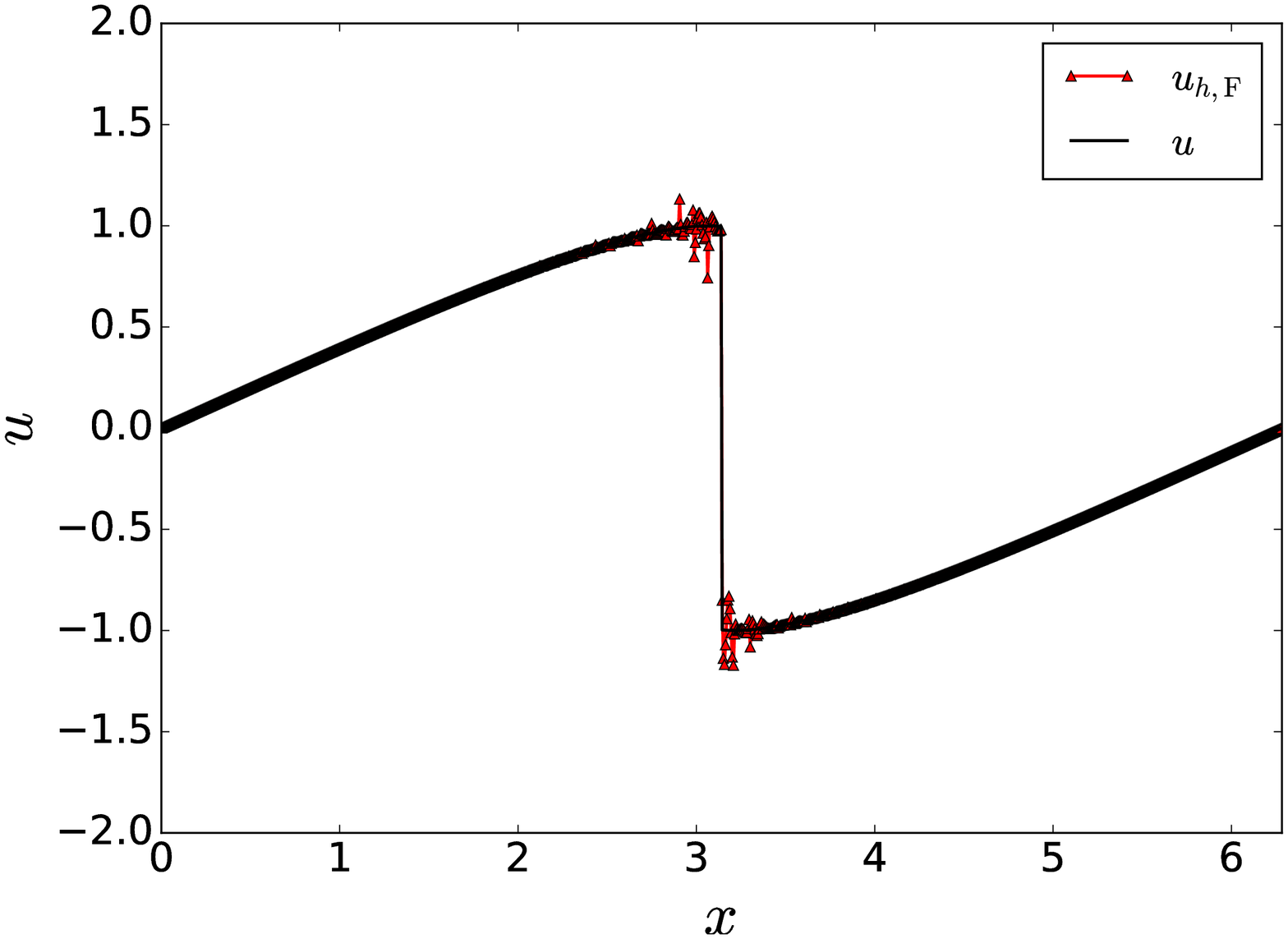}
		\caption{$\nu = -\frac{1}{4}$.}	
	\end{subfigure}
	\caption{Solution profiles at $T = 1.5$, $p = 5$, $ k = 4$, $N = 80$.}
\end{figure}

\section{Conclusion}\label{sec-concl}
\setcounter{equation}{0}
\setcounter{figure}{0}
\setcounter{table}{0}

In this paper, we use the idea of superviscosity from spectral methods to enforce strong stability of explicit RK methods. The main efforts are put on linear ODE systems. The superviscosity is imposed either by modifying the semi-negative operator (the modified method), or by post-processing the solution (the filtering method). For both methods, we estimate the energy difference between the stabilized solution and the original solution. Then it is shown that by introducing certain amount of superviscosity, the RK method can achieve strong stability while maintaining the original order of accuracy. For diffusive superviscosity, we derive a sharp bound on the minimum superviscosity that is needed to ensure strong stability. For the general dispersive-diffusive superviscosity,  the bound we provide is verified to be sharp numerically. A filtering method is then presented for nonlinear problems. Finally numerical tests are given to examine the sharpness on our estimates, the accuracy, the effect on energy evolution of the modified and filtering methods.

\end{document}